\documentclass[11pt]{article}
\usepackage{definitions}

\begin{document}

\title{An $O(s^r)$-Resolution ODE Framework for Understanding Discrete-Time Algorithms and Applications to the Linear Convergence of Minimax Problems}

\author{Haihao Lu\thanks{The University of Chicago Booth School of Business
({mailto:  haihao.lu@uchicagobooth.com}).}
}
\date{{\scriptsize (final version May 2021, first version Jan 2020)}}


\maketitle

\begin{abstract}
\redd{
There has been a long history of using ordinary differential equations (ODEs) to understand the dynamics of discrete-time algorithms (DTAs). Surprisingly, there are still two fundamental and unanswered questions:
(i) it is unclear how to obtain a \emph{suitable} ODE from a given DTA, and (ii) it is unclear the connection between the convergence of a DTA and its corresponding ODEs. In this paper, we propose a new machinery -- an $O(s^r)$-resolution ODE framework -- for analyzing the behavior of a generic DTA, which (partially) answers the above two questions. The framework contains three steps: 1. To obtain a suitable ODE from a given DTA, we define a hierarchy of $O(s^r)$-resolution ODEs of a DTA parameterized by the degree $r$, where $s$ is the step-size of the DTA. We present a principal approach to construct the unique $O(s^r)$-resolution ODEs from a DTA; 2. To analyze the resulting ODE, we propose the $O(s^r)$-linear-convergence condition of a DTA with respect to an energy function, under which the $O(s^r)$-resolution ODE converges linearly to an optimal solution; 3. To bridge the convergence properties of a DTA and its corresponding ODEs, we define the properness of an energy function and show that the linear convergence of the $O(s^r)$-resolution ODE with respect to a proper energy function can automatically guarantee the linear convergence of the DTA.
}


To better illustrate this machinery, we utilize it to study three classic algorithms -- gradient descent ascent (GDA), proximal point method (PPM) and extra-gradient method (EGM) -- for solving the unconstrained minimax problem $\min_{x\in\RR^n} \max_{y\in \RR^m} L(x,y)$. Their $O(s)$-resolution ODEs explain the puzzling convergent/divergent behaviors of GDA, PPM and EGM when $L(x,y)$ is a bilinear function, and showcase that the interaction terms help the convergence of PPM/EGM but hurts the convergence of GDA. Furthermore, their $O(s)$-linear-convergence conditions not only unify the known scenarios when PPM and EGM have linear convergence, but also showcase that these two algorithms exhibit linear convergence in much broader contexts, including when solving a class of nonconvex-nonconcave minimax problems. Finally, we show how this ODE framework can help design new optimization algorithms for minimax problems, by studying the difference between the $O(s)$-resolution ODE of GDA and that of PPM/EGM.

\end{abstract}


\section{Introduction}

\begin{figure}[!tbp]
\centering
  \begin{subfigure}[b]{0.48\textwidth}
    \includegraphics[width=\textwidth, trim=0 0 0 0, clip]{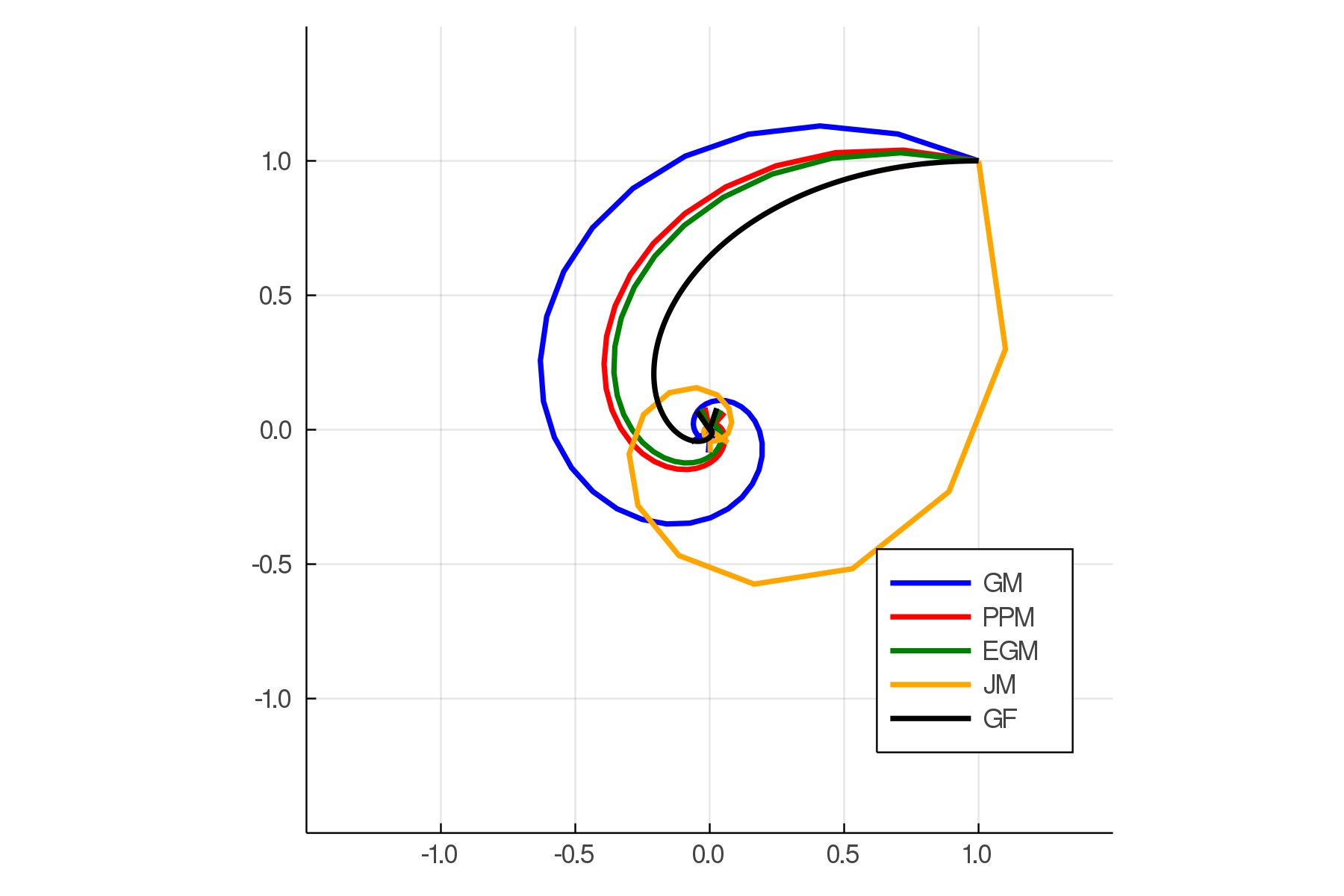}
    \caption{The trajectories of different algorithms for solving $\min_x \max_y \frac{1}{2} x^2 + 2xy -\frac{1}{2} y^2$ with step-size $s=0.1$ and initial solution $(1,1)$.}
  \end{subfigure}
  \hspace{0.2cm}
  \begin{subfigure}[b]{0.48\textwidth}
    \includegraphics[width=\textwidth, trim=0 0 0 0, clip]{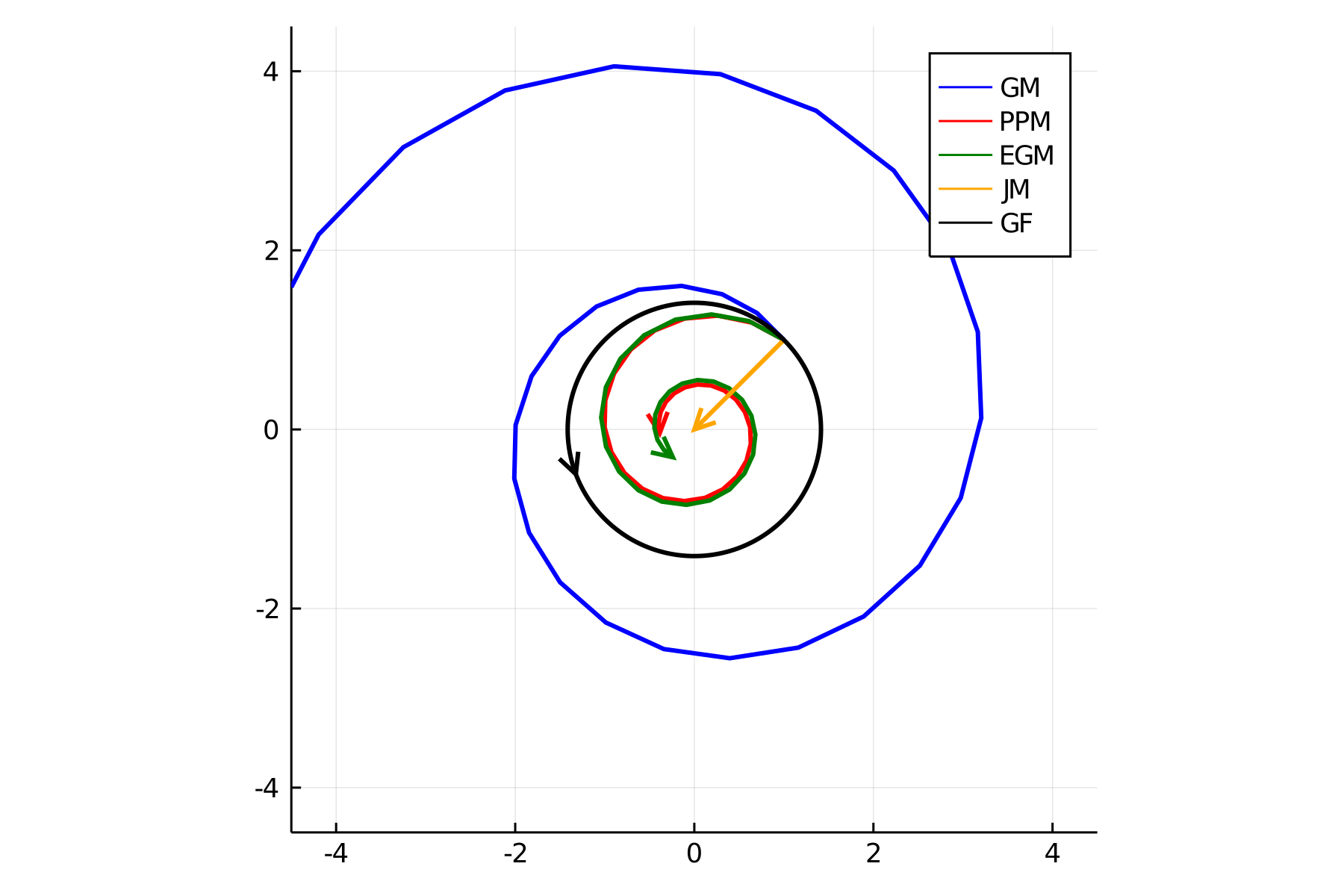}
    \caption{The trajectories of different algorithms for solving $\min_x \max_y xy $ with step-size $s=0.3$ and initial solution $(1,1)$.}
  \end{subfigure}
  \caption{Illustration of the behaviors of GDA, PPM, EGM, JM  (Jacobian method introduced later in \eqref{eq:new_al}) and GF for solving minimax problems in the two scenarios when $L(x,y)$ is strongly convex and when $L(x,y)$ is bilinear.}
  \label{fig:discrete}
\end{figure}

There has been a long history of using ordinary differential equations (ODEs) to understand the dynamics of discrete-time algorithms (DTAs) \cite{helmke2012optimization,schropp2000dynamical,fiori2005quasi}. Recently, the seminal work \cite{su2016differential} triggered a renewed spark on this line of research. The ODE perspective to understand DTAs has two major advantages: the convergence analysis for ODEs is usually more straight-forward than that for DTAs; and the advanced analytical tools from ODE literature can help provide more fundamental intuitions on the behaviors of DTAs~\cite{su2016differential}. However, there are still two fundamental unanswered questions when utilizing this approach:

\begin{itemize}
    \item How to obtain a suitable ODE from a given DTA? Indeed, there can be multiple ODEs that correspond to the same DTA, depending on how to take the continuous limit~\cite{shi2018understanding}. While the easiest approach to construct an ODE from a DTA is by simply letting the step-size $s$ go to $0$, the resulting ODEs may not be able to distinguish different DTAs, and even worse, the trajectories of the DTA and such ODEs can be topologically different with any positive step-size $s$ (see for example Figure \ref{fig:discrete} (b)).
    \item What is the connection between the convergence of a DTA and the convergence of its corresponding ODE? Although the convergence analysis for ODEs, in many cases, is straight-forward, translating it back to the convergence of DTAs (if it is possible) can be highly nontrivial.
\end{itemize}


For example, the derivation of the ODE corresponding to Nesterov's accelerated method in~\cite{su2016differential, shi2018understanding} is somewhat ``informal'', and requires some good mathematical intuitions on how and where to perform the Taylor expansion; in the meantime, the convergence guarantees of the DTAs require an independent and highly technical analysis on top of analysis for the corresponding ODEs~\cite{shi2018understanding}.

In this paper, we propose an $O(s^r)$-resolution ODE framework to analyze the behavior of DTAs, which (partially) resolves the above two questions. We study a generic DTA with iterate update:
\begin{equation}\label{eq:discrete-update}
    z^+ = g(z, s) \ ,
\end{equation}
where $z$ is the iterate input, $z^+$ is the iterate output,  $s$ is the step-size of the algorithm, and $g(z,s)$ is a sufficiently smooth
function in $z$ and $s$, which satisfies that $g(z,0)=z$ (i.e. the current solution does not move if the step-size $s=0$). We propose an $O(s^r)$-resolution ODE framework for analyzing a DTA \eqref{eq:discrete-update}, which contains the following three key steps:
\redd{
\begin{enumerate}
    \item \textbf{Obtain an ODE from a DTA}: Choose a suitable degree $r$, and perform the $r$-th degree ODE expansion of the DTA to obtain its $O(s^r)$-resolution ODE (see Section \ref{sec:high-resolution}). The value of $r$ should be chosen so that the $O(s^r)$-resolution ODE is capable to characterize the major (convergent) behaviors of the DTA.
    \item \textbf{Analyze the ODE}: Choose an energy function, and obtain the $O(s^r)$-linear-convergence conditions of the DTA, under which the resulting $O(s^r)$-resolution ODE linearly converges to an optimal solution with respect to this energy function (see Section \ref{sec:linear-convergence-condition}). 
    \item \textbf{Translate the convergent results back to the DTA}: Under mild conditions, the $O(s^r)$-linear-convergence conditions obtained in the previous step can automatically guarantee the convergence of the DTA if the energy function chosen in the previous step is \emph{proper} (see Section \ref{sec:ODEtoDTA}), and it can also motivate a direct convergence analysis in the discrete-time space (see Section \ref{sec:convergence}). These connections between the DTA and the ODEs heavily rely on the construction of the $O(s^r)$-resolution ODE.
\end{enumerate}}

This framework is inspired by the recent work of the high-resolution ODE for analyzing the difference between Nesterov's accelerated method and heavy ball method \cite{shi2018understanding}. The key differences between our framework and that in \cite{shi2018understanding} are: (i) we propose the $r$-th degree ODE expansion of a DTA to obtain its corresponding $O(s^r)$-resolution ODE, while their informal derivation of the $O(s)$-resolution ODE of momentum methods in \cite{shi2018understanding} may not be easily generalized to other algorithms or to higher order resolution ODEs; (ii) we fix the energy function first and then study for what class of problems the ODE has linear convergence with respect to this energy function, while they focus on constructing a decaying energy function under the standard convexity conditions; (iii) under mild conditions, the linear convergence of the $O(s^r)$-resolution ODEs can automatically guarantee the linear convergence of the DTA, while their analysis of the DTA is independent of the ODE analysis and it can be highly non-trivial.

To further illustrate the ideas of the $O(s^r)$-resolution ODE framework, we study the following unconstrained minimax problem as an example:
\begin{equation}\label{eq:POI}
\min_{x\in\RR^n} \max_{y\in \RR^m} L(x,y) \ ,    
\end{equation}
where $L(x,y)\in \RR^m \times \RR^n \rightarrow \RR$ is a sufficiently differentiable function. The goal is to design first-order methods to find a stationary point (equivalently a first-order Nash equilibrium) $(x^*,y^*)$ of \eqref{eq:POI} such that
\begin{equation}
\nabla_x L(x^*,y^*)=0 \text{ and } \nabla_y L(x^*,y^*)=0 \ .
\end{equation}
Define $z=(x,y)\in \RR ^{n+m}$ and $F(z)= [\nabla_x L(x,y), -\nabla_y L(x,y)]\in \RR ^{n+m}$, then $z^*=(x^*,y^*)$ is a stationary point of \eqref{eq:POI} iff $F(z^*)=0$. We will utilize $z$ and $F(z)$ throughout the paper for notational convenience.

Minimax problem \eqref{eq:POI} has many applications, including but not limited to: generative adversarial networks~\cite{goodfellow2014generative}, robust optimization~\cite{ben2009robust,bertsimas2011theory}, Lagrangian formulation of constrained convex optimization~\cite{Rock}, supervised learning~\cite{zhang2017stochastic}, matrix factorization~\cite{bach2008convex}, PID robust control~\cite{hast2013pid}, etc.

Here we study the following three classic algorithms for solving \eqref{eq:POI}, and focus on their linear convergence rate: 
\begin{itemize}
    \item Gradient Descent Ascent (GDA): 
    \begin{equation}\label{eq:gd}
    z_{+}=z - s F(z)\ ,
    \end{equation}
    \item Proximal Point Method (PPM):
    \begin{equation}\label{eq:ppm}
    z_{+}=z - s F(z_{+})\ ,    
    \end{equation}
    \item Extra-Gradient Method (EGM) (it is a special case of Mirror Prox Algorithm \cite{nemirovski2004prox}): 
    \begin{equation}\label{eq:mp}
    \tz=z - s F(z), z_{+}=z - s F(\tz)\ ,
    \end{equation}
\end{itemize}
where $s$ is the step-size of each algorithm.

There have been extensive studies on analyzing the computational guarantees of the above three algorithms for solving \eqref{eq:POI}. Essentially, previous works show that linear convergence occurs under one of the following two scenarios:
\begin{itemize}
    \item[] (i) $L(x,y)$ is strongly convex-strongly concave, i.e. $L(x,y)$ is strongly convex in $x$ and strongly concave in $y$;
    \item[] (ii) $L(x,y)=x^T B y$ is a bilinear function. 
\end{itemize}

More specifically, it has been shown that all three algorithms have linear convergence in Scenario  (i), but there is a puzzling phenomenon in Scenario (ii): while PPM and EGM converge linearly, GDA diverges ~\cite{bauschke2011convex,du2019linear,rockafellar1976monotone,tseng1995linear, wang2017exploiting,liang2019interaction}. See Figure \ref{fig:discrete} for examples of the above behaviors. A more detailed literature review is presented in Section \ref{sec:literature}.


Indeed, GDA, PPM and EGM are highly related. When the step-size $s$ goes to $0$, one can show that all of these three algorithms result in the same continuous-time system --- gradient flow (GF),
\begin{equation}\label{eq:GF}
    \dZ= - F(Z) \ .
\end{equation}
Moreover, they all share similar trajectories towards a stationary point of \eqref{eq:POI} in Scenario (i) (See Figure \ref{fig:discrete} (a) for an example). However, it is a mystery to see that these three algorithms exhibit topologically different behaviors in Scenario (ii) -- GDA diverges, PPM and EGM converges to a stationary point of \eqref{eq:POI}, and GF keeps oscillating and never converge nor diverge (see Figure \ref{fig:discrete} (b) for an example). This work provides an intuitive explanation of the above puzzling behaviors via the $O(s)$-resolution ODEs of GDA, PPM and EGM. As we will see later, such strange behaviors are due to a multi-scale phenomenon: The linear convergence in Scenario (i) is an $O(1)$-scale behavior; the three methods result in the same $O(1)$-resolution ODE system (i.e., GF), thus they share similar convergent behaviors, following the path of GF. On the other hand, Scenario (ii) is a limiting case when an $O(s)$-perturbation of the dynamic can dramatically change the behavior of GF, thus we need to look at $O(s)$-resolution approximation of the discrete-time algorithms in order to understand their trajectories. As we will show in Section \ref{sec:high-resolution}, the $O(s)$-resolution ODEs of GDA, PPM and EGM contain an extra term -- $\frac{s}{2}\nabla F(Z) F(Z)$ with different signs on top of the dynamics of GF, which is the fundamental reason of the above convergent/divergent behaviors of the GDA and PPM/EGM. Furthermore, while both PPM and EGM share similar trajectories in Scenario (ii) (since they share the same $O(s)$-resolution ODE), they have subtle frequency discrepancy. This is an $O(s^2)$-behavior, which can be explained by the difference in their $O(s^2)$-resolution ODEs. Motivated by the difference between the $O(s)$-resolution ODEs of GDA and and that of PPM/EGM, we design a new algorithm, Jacobian method (JM), for minimax problems, which can avoid spiral and go directly to the minimax solution when the objective $L(x,y)$ is bilinear.

Furthermore, the above two scenarios when PPM/EGM has linear convergence are disconnected, in particular, compared with the clean and unified linear convergence results in convex optimization literature~\cite{nesterovBook}. Recall that in the classic convex optimization theory, gradient-based methods with a reasonably small step-size $s$ find a solution within $\varepsilon$ optimality gap in $O(\frac{1}{s\mu}\log\frac{1}{\varepsilon})$ iterations, where $\mu$ is the strong convexity constant of the objective function defined by the Hessian of the objective function~\cite{nesterovBook}. However, to the best of our knowledge, there is a lack of such a simple constant which naturally characterizes the linear convergence rate of different algorithms for solving minimax problem \eqref{eq:POI}.
Here, the $O(s)$-resolution ODEs of PPM and EGM inspire us to introduce the $O(s)$-linear-convergence constant $\rho(s)$, which is defined by the Hessian of $L(x,y)$ and the step-size $s$ of the algorithm, and similar to the classic convex optimization, PPM and EGM find a solution $z$ such that $\|F(z)\|^2\le\varepsilon$ in $O(\frac{1}{s\rho(s)}\log\frac{1}{\varepsilon})$ iterations with a reasonably small step-size $s$. This constant $\rho(s)$ not only unifies the known linear convergence rate of PPM and EGM in the above two classic scenarios, but also showcases that these two algorithms exhibit linear convergence in broader contexts, including a class of nonconvex-nonconcave minimax problems  (see Example \ref{ex:3}-\ref{ex:6} in Section \ref{sec:linear-convergence-condition}). Indeed, such analysis clearly shows that the interaction term in $L(x,y)$ helps the convergence of PPM and EGM, but hurts the convergence of GDA.

In the rest of this section, we present the related literature and a summary of the contributions of this work.

\subsection{Related Literature} \label{sec:literature}
In the seminal work~\cite{rockafellar1976monotone}, Rockafellar studied PPM for solving monotone variational inequalities. For minimax problems \eqref{eq:POI} (as a special case of variational inequality), his results imply that PPM has local linear convergence under the conditions that (a) the solution to \eqref{eq:POI} is unique, (b) the function $F:\RRmn\rightarrow\RRmn$ is invertible around $0$, and (c) $F^{-1}$ is Lipschitz continuous around $0$, which are satisfied in Scenario (i). Moreover, \cite{rockafellar1976monotone} further shows that PPM has global linear convergence in Scenario (i). Later on, Tseng~\cite{tseng1995linear} shows that both PPM and EGM have a linear convergence rate for solving variational inequality when certain complicated conditions are satisfied, and these conditions are satisfied for solving the minimax problem \eqref{eq:POI} in Scenario (i) and in Scenario (ii). In 2004, Nemirovski~\cite{nemirovski2004prox} proposes Mirror Prox algorithm (a special selection of the prox function recovers EGM), which first shows that EGM has $O(\frac{1}{\varepsilon})$ sub-linear convergence rate for solving convex-concave minimax problems over a compact set.

There are several works that study the special case of \eqref{eq:POI} when the minimax function has bilinear interaction terms, i.e., $L(x,y)=f(x)+x^T B y - g(y)$ where $f(\cdot)$ and $g(\cdot)$ are both convex functions. The most influential algorithms for solving the above bilinear interaction minimax problems are perhaps Nesterov's smoothing~\cite{nesterov2005smooth}, Monteiro's hybrid proximal extragradient method~\cite{monteiro2010complexity}, Douglas-Rachford splitting (a special case is Alternating Direction Method of Multipliers (ADMM))~\cite{douglas1956numerical,eckstein1992douglas} and Primal-Dual Hybrid Gradient Method (PDHG)~\cite{chambolle2011first} (the last two are recently shown to be equivalent under preconditioning~\cite{o2018equivalence}). Moreover, ADMM and PDHG also have linear convergence under different types of conditions, but a major difference between these two algorithms and the methods studied in this paper is that these two algorithms do the primal update and the dual update sequentially, while PM, PPM and EGM do the primal update and the dual update simultaneously.

More recently, minimax problems have gained attention in machine learning community, perhaps mainly due to the study on Generative Adversarial Networks (GANs). \cite{daskalakis2018training} studies the Optimistic Gradient Descent Ascent (OGDA) designing for training GANs, and shows that OGDA converges linearly for bilinear minimax problems with additional assumptions that the matrix $B$ is square and full rank (it is thus a special case of Scenario (ii)). \cite{mokhtari2020unified} shows that OGDA, EGM both approximate PPM (indeed, EGM is an approximation to PPM was first shown in Nemirovski's earlier work \cite{nemirovski2004prox}), and further showed that these three algorithms have a linear convergence rate when $L(x,y)$ is strongly convex-strongly concave (Scenario (i)) or when $L(x,y)$ is bilinear with square and full rank matrix $B$ (again, a special case of Scenario (ii)). See~\cite{mokhtari2020unified} for a more detailed literature review on recent results on OGDA. Although we do not study OGDA in this paper, we do not see any reason that the techniques and results developed herein cannot be used to analyze the performance of OGDA or other types of inexact PPM. 

Another recent line of research on continuous optimization is to understand the optimization methods from the continuous-time dynamical system perspective. Su, Boyd and Candes~\cite{su2016differential} presents the $O(1)$-resolution ODE system of Nesterov's accelerated method~\cite{nesterov1983method} for convex optimization, which provides a new explanation of why Nesterov's method can speed up the convergence rate of gradient-based methods. Later on, Lagrangian and Hamiltonian frameworks are proposed to understand the acceleration phenomenon and generate a larger class of accelerated methods~\cite{wilson2016lyapunov, wibisono2016variational}. More recently, \cite{shi2018understanding} proposes an $O(s)$-resolution ODE system that explains the different behaviors between Nesterov's accelerated method and heavy-ball method, even though both algorithms share the same $O(1)$-resolution ODE. Refer to~\cite{shi2018understanding} for a more detailed literature review on this line of research. 

Lastly, we want to mention that the multi-scale expansion of the perturbation of a continuous-time ODE system has been well studied in physics and in applied mathematics~\cite{pedlosky2013geophysical,weinan2011principles}. 
\redd{
\subsection{Summary of Contributions}
We present a new machinery -- an $O(s^r)$-resolution ODE framework -- for analyzing the behavior of a generic discrete-time algorithm, and apply it to unconstrained minimax problems:
\begin{itemize}
    \item From DTAs to ODEs: Given a DTA, we introduce its $O(s^r)$-resolution ODE (Definition \ref{def:Osr-ODE}), and propose an $r$-th degree ODE expansion to obtain the unique $O(s^r)$-resolution ODE (Theorem \ref{thm:high-resolution}).
    \item Analyze the ODEs: We propose to study $O(s^r)$-linear-convergence conditions with respect to an energy function under which the $O(s^r)$-resolution ODE converges linearly (Definition \ref{def:linear-conv-condition}).
    \item From ODEs to DTAs: We introduce the properness of an energy function to study the $O(s^r)$-resolution ODE of a DTA (Definition \ref{def:proper_energy}), and present easy-to-check sufficient conditions (Theorem \ref{thm:gradient-based-ODE}).
    We show that with a \emph{proper} choice of the energy function, the linear convergence of the $O(s^r)$-resolution ODE can automatically guarantee that the DTA has the same linear convergence rate (Theorem \ref{thm:gradient-based-ODE-bound}).
   
    \item We utilize the above framework to study GDA, PPM and EGM for solving minimax problem \eqref{eq:POI}. When $L(x,y)$ is a bilinear function, the closed-form solutions to their $O(s)$-resolution ODEs explain the puzzling behaviors of the three algorithms. Furthermore, the closed-form solutions to the $O(s^2)$-resolution ODEs of PPM and EGM explain their subtle frequency discrepancy (Section \ref{sec:understanding_behaviors}).
    \item We propose to study the energy function $\frac{1}{2}\|F(z)\|^2$ for analyzing the convergence of PPM and EGM for minimax problems, and we show $\frac{1}{2}\|F(z)\|^2$ is a proper energy function. Using the above framework, we introduce the $O(s)$-linear-convergence condition of PPM and EGM for solving \eqref{eq:POI}, which not only unifies the linear convergence results in previous works, but also showcases that PPM and EGM exhibit linear convergence in broader contexts (Section \ref{sec:ODEtoDTA} and Section \ref{sec:convergence}). 
    \item Inspired by the difference between the $O(s)$-resolution ODE of PPM/EGM and that of GDA, we introduce a new algorithm, Jacobian Method (JM), which avoids the spiral and can go directly towards the stationary point for minimax problems with sufficient interaction terms (Section \ref{sec:design-new-al}).

\end{itemize}
}










\subsection{Notations} We use $\ell_2$-norm throughout the paper, namely, $\|c\|=\sqrt{\sum_i c_i^2}$ for any vector $c$, and $\|M\|=\max_{x,y} \frac{y^T M x}{\|x\|\|y\|}$ for any matrix $M$. For a symmetric matrix $M$, $\lambda_{\min}(M)$ is the minimal eigenvalue of $M$. For a positive-semidefinite matrix $M$, $\lambda_{\min}^+(M)$ is the minimal non-zero eigenvalue of $M$. We denote $A(z)=\nabla_{xx}L(x,y)$, $B(z)=\nabla_{xy}L(x,y)$, $C(z)=-\nabla_{yy}L(x,y)$, then $\nabla F(z)=\left[\begin{smallmatrix} A(z) & B(z) \\ -B(z)^T & C(z) \end{smallmatrix}\right]$. We also use $A, B, C$ to represent $A(z), B(z), C(z)$ if they do not cause any misunderstandings. $\text{Conv}(S)$ refer to the convex hull of a set $S$.





\section{From DTAs to ODEs: The $O(s^r)$-Resolution ODE of a DTA}\label{sec:high-resolution}
In this section, we introduce the $r$-th degree ODE expansion of a DTA to obtain the unique $O(s^r)$-resolution ODE of the DTA. Based on that, we obtain the $O(1)$-resolution ODEs of GDA, PPM and EGM, which explains the convergent behaviors of these three algorithms in Scenario (i); we obtain the $O(s)$-resolution ODEs of GDA, PPM and EGM, whose solutions explain the puzzling divergent/convergent behaviors of the three algorithms in Scenario (ii); and we obtain the $O(s^2)$-resolution ODEs of PPM and EGM, whose solutions explain their frequency discrepancy  in Scenario (ii). Finally, we discuss how the $O(s)$-resolution ODE framework can help design new algorithms.



\vspace{0.2cm}

\subsection{The $O(s^r)$-Resolution ODE}\label{sec:Os-resolution-ODE}
First, let us formally define an $O(s^r)$-resolution ODE of a DTA:

\begin{mydef}\label{def:Osr-ODE}
We say an ODE system with the following format
\begin{equation}\label{eq:rthode}
    \dZ = f^{(r)}(Z,s) := f_{0}(Z)+s f_{1}(Z) + \cdots + s^r f_{r}(Z)
\end{equation}
the $O(s^r)$-resolution ODE of the DTA with iterate update \eqref{eq:discrete-update} if it satisfies that for any $z$ and $z^+=g(z,s)$,
\begin{equation}\label{eq:small_term}
    \|Z(s)-z^+\| = o(s^{r+1}) \ , \footnote{Recall that the $o$ notation in Equation \eqref{eq:small_term} means $\lim_{s\rightarrow 0} \frac{\|Z(s)-z^+\|}{s^{r+1}}=0$.}
\end{equation}
where $Z(s)$ is the solution obtained at $t=s$ following the ODE \eqref{eq:rthode} with initial solution $Z(0)=z$.
\end{mydef}

\vspace{0.2cm}

Next, we describe how to obtain the $O(s^r)$-resolution ODE from the discrete-time update function $g(z,s)$, and we call this process the \underline{$r$-th degree ODE expansion of a DTA}. Before that, let us introduce some new notations:

Suppose the function $g(z,s)$ is $(r+1)$-th order differentiable over $s$ for any $z$, then by Taylor expansion of $g(z,s)$ over $s$ at $s=0$, we obtain
\begin{equation}\label{eq:expan_g}
    g(z,s)=\sum_{j=0}^{r+1} \frac{1}{j!} \left.\frac{\partial^j g(z, s)}{\partial s^j}\right|_{s=0} s^j + o(s^{r+1})=\sum_{j=0}^{r+1} g_j(z) s^j + o(s^{r+1}) \ ,
\end{equation}
where $g_j(z):=\frac{1}{j!} \left.\frac{\partial^j g(z, s)}{\partial s^j}\right|_{s=0}$ is the $j$-th coefficient function in the above Taylor expansion.

Suppose $f_i(Z)$ in \eqref{eq:rthode} is $(r+1)$-th order differentiable for $i=0,\ldots,r$, then $\frac{d^j}{dt^j} Z$ exists for any $j=0,\ldots,r+1$, and it is a $jr$-th order polynomial in $s$. Let us define $h_{j,i}(Z)$ as the coefficient function of $s^i$ in the expansion of $\frac{d^j}{dt^j} Z$, i.e., 
\begin{equation}\label{eq:djdt}
\frac{d^j}{dt^j} Z=\sum_{i=0}^{r+1}  h_{j,i}(Z) s^i + o(s^{r+1}).
\end{equation}
Substituting \eqref{eq:rthode} into \eqref{eq:djdt} and comparing the coefficient function of $s^0, s^1, \ldots, s^i$ on both sides of \eqref{eq:djdt}, we have that $h_{j,i}(Z)$ is a function of $f_0(Z),\ldots, f_i(Z)$ for any $0\le i\le r,\ 0\le j\le r+1$. Moreover, it holds that
\begin{itemize}
    \item when $j=0$, we have $\frac{d^0}{dt^0} Z = Z$, thus $h_{0,0}(Z)=Z$ and $h_{0,i}(Z)=0$ for $i= 1,2,\ldots,r$;
    \item when $j=1$, we have $\frac{d^1}{dt^1} Z = f^{(r)}(Z,s)$, thus $h_{1,i}(Z)=f_i(Z)$ for $i=0,\ldots, r$;
    \item when $j=2$, we have $\frac{d^2}{dt^2} Z = \nabla_z f^{(r)}(Z,s) f^{(r)}(Z,s)$, thus $h_{2,i}(Z)=\sum_{l=0}^i \nabla f_{i-l}(Z) f_l   (Z)$ for $i=0,\ldots, r$;
    \item more generally, the functions $h_{j,i}(Z)$ can be computed recursively by taking the derivative over $t$ in \eqref{eq:djdt} and comparing the corresponding terms as
\begin{equation}\label{eq:recursive_compute_h}
    h_{j+1,i}(Z)=\sum_{l=0}^i \nabla h_{j,l}(Z) h_{1, i-l}(Z) \ .
\end{equation}
\end{itemize}

\vspace{0.2cm}
The next theorem presents the $r$-th degree ODE expansion of a DTA, through which we obtain its corresponding $O(s^r)$-resolution ODE:
\begin{thm} 
\label{thm:high-resolution}
Consider a DTA with iterate update $z_{+}=g(z,s)$, where $g(z,0)=z$ and $g(z,s)$ is sufficiently differentiable in $s$ and in $z$. Then its $O(s^r)$-resolution ODE is unique, and the $i$-th coefficient function in the $O(s^r)$-resolution ODE can be obtained recursively by
\begin{equation}\label{eq:obtainf}
f_{i}(Z)=g_{i+1}(Z)-\sum_{l=2}^{i+1}  \frac{1}{l!} h_{l,i+1-l}(Z) \ , \emph{ for } i=0,1,\ldots, r,
\end{equation}
where $h_{l,i+1-l}(Z)$ is defined in \eqref{eq:djdt} and it is a function of $f_{0}(Z),\ldots,f_{i-1}(Z)$ for $2\le l\le i+1$. 
\end{thm}

\textbf{Proof. } Suppose there exists an $O(s^r)$-resolution ODE \eqref{eq:rthode} of the DTA with iterate update $z^+=g(z,s)$. By Taylor expansion of $Z(t)$ at $t=0$, we obtain that
\begin{align}\label{eq:expan_Zs}
    \begin{split}
        Z(s) &=  \sum_{j=0}^{r+1} \frac{1}{j!} \frac{d^j}{dt^j} Z(0)s^j + o(s^{r+1}) \\
        & = \sum_{j=0}^{r+1}  \frac{1}{j!} s^j \sum_{i=0}^{r+1} h_{j,i} (Z(0)) s^i + o(s^{r+1}) \\
        & = \sum_{j=0}^{r+1} \sum_{l=0}^{j}  \frac{1}{l!} h_{l,j-l}(Z(0)) s^j + o(s^{r+1}) \ , \\
        & = \sum_{j=0}^{r+1} \sum_{l=0}^{j}  \frac{1}{l!} h_{l,j-l}(z) s^j + o(s^{r+1}) \ , 
    \end{split}
\end{align}
where the second equality uses \eqref{eq:djdt} and the last equality is from $Z(0)=z$. Notice that the $O(s^r)$-resolution ODE satisfies \eqref{eq:small_term}, thus the coefficient functions of $s^j$ in the expansion \eqref{eq:expan_g} and in the expansion \eqref{eq:expan_Zs} must be the same. Therefore, it holds for $0\le j\le r+1$ that
\begin{equation}\label{eq:kth_tempo}
    \sum_{l=0}^{j}  \frac{1}{l!} h_{l,j-l}(z) =g_j(z) \ .
\end{equation}
By rearranging \eqref{eq:kth_tempo} and noticing $h_{0,j+1}=0$ and $h_{1,j}(z)=f_{j}(z)$, we have for any $1\le j\le r$ that
\begin{equation}\label{eq:recursive}
    f_{j}(z)=h_{1,j}(z)=g_{j+1}(z)-\sum_{l=2}^{j+1}  \frac{1}{l!} h_{l,j+1-l}(z) \ ,
\end{equation}
In particular, when $j=0$ we have that $f_0(z)=h_{1,0}(z)=g_1(z)-h_{0,1}(z)=g_1(z)$. Notice that $h_{l,j+1-l}(z)$ is a function of $f_{0}(z), f_{1}(z),\ldots, f_{j-1}(z)$ for any $2\le l \le j+1$, thus the right-hand side of \eqref{eq:recursive} is a function of $g_{j+1}(z), f_{0}(z), f_{1}(z),\ldots, f_{j-1}(z)$, which provides a recursive way to define $f_j(z)$ from $g_1(z),\ldots,g_{j+1}(z)$.

The above process also guarantees that the obtained ODE \eqref{eq:rthode} with coefficient function $f_j(z)$ from \eqref{eq:recursive} satisfies \eqref{eq:small_term}, thus it is indeed an $O(s)$-resolution ODE of the DTA \eqref{eq:discrete-update}. Furthermore these $f_j(z)$ is uniquely defined by $g_1(z),\ldots,g_{j+1}(z)$ through \eqref{eq:recursive}, thus the $O(s^r)$-resolution ODE of a DTA is unique.\qed


\begin{rem}\label{rem:high-order-difference}
Indeed,
the $O(s)$-resolution ODE results in a stronger bound when $g(z,s)$ is sufficiently smooth:
\begin{equation}\label{eq:small_term_strong}
    \|Z(s)-z^+\| = O(s^{r+2}) \ . \footnote{Recall that the $O(\cdot)$ notation in Equation \eqref{eq:small_term_strong} is equivalent to that there exists a constant $C$ such that $\lim_{s\rightarrow 0} \frac{\|Z(s)-z^+\|}{s^{r+2}}\le C$.}
\end{equation}
This can be simply obtained from the proof of Theorem \ref{thm:high-resolution} by replacing $o(s^{r+1})$ to $O(s^{r+2})$. 
\end{rem}

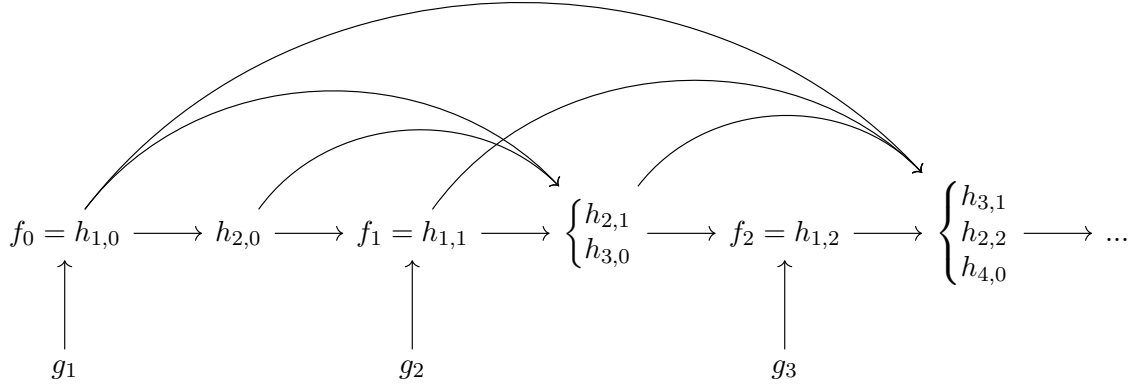
\begin{figure}
    \centering
     \begin{tikzcd}
  f_0=h_{1,0} \arrow[r] \arrow[rrr, bend left=50, ""{name=U, below}] \arrow[rrrrr, bend left=50, ""{name=U, below}] &  h_{2, 0} \arrow{r} \arrow[rr, bend left=50, ""{name=U, below}] & f_1=h_{1,1} \arrow[r] \arrow[rrr, bend left=50, ""{name=U, below}] & \left\{\begin{matrix}
    h_{2,1} \\ h_{3,0} 
\end{matrix} \right. \arrow[r] \arrow[rr, bend left=50, ""{name=U, below}] & f_2=h_{1,2} \arrow[r] & \left\{\begin{matrix}
    h_{3,1} \\ h_{2,2} \\ h_{4,0} 
\end{matrix} \right. \arrow[r] & ... \\
g_1 \arrow{u}&& g_2 \arrow{u} &  & g_3 \arrow{u}& 
\end{tikzcd}
    \caption{\red{The logic flow of computing the high-resolution ODE (i.e., the coefficient functions $f_j=h_{1,j}$ in \eqref{eq:rthode}) recursively from the DTA update (i.e., $g_1, g_2, ...$).}}
    \label{fig:logic_flow}
\end{figure}


\red{
Figure \ref{fig:logic_flow} plots the logic flow to compute the $O(s^r)$-resolution ODE recursively from the Taylor coefficient functions $\{g_j\}$ of a DTA. Suppose we know $h_{i,j}$ for $i+j\le k$. Then, we can compute $h_{i,j}$ for $i+j = k+1$ as follows: We obtain $h_{i,j}$ for $i+j=k+1$ and $i\ge 2$ using \eqref{eq:recursive_compute_h}, and then we obtain $f_k=h_{1, k}$ using \eqref{eq:recursive}.

}

Following Theorem \ref{thm:high-resolution}, we present a conjecture:
\begin{con}\label{con:inf}
Under certain regularity conditions on $g(z,s)$ and $s$ (for example, $g(z,s)$ is infinitely differentiable, $s$ needs to be reasonably small, etc), the infinite sum in the right-hand-side of
\begin{equation*}
    f^{(\infty)}(Z,s):= \sum_{i=0}^{\infty} f_i(Z)s^i
\end{equation*}
converges for any $Z$, where $f_i(Z)$ is defined recursively by \eqref{eq:obtainf}. Furthermore, for any $z$ and $z^+=g(z,s)$, it holds that
\begin{equation*}\label{eq:small_term_inf}
    Z(s)=z^+ \ ,
\end{equation*}
where $Z(s)$ is the solution obtained at $t=s$ following from the ODE system
\begin{equation}\label{eq:infode}
    \dZ=f^{(\infty)}(Z,s)
\end{equation}
with initial solution $Z(0)=z$. \qed
\end{con}
Suppose Conjecture \ref{con:inf} holds, then the ODE system \eqref{eq:infode} can fully characterize the DTA with iterate update \eqref{eq:discrete-update}. In particular, suppose $z_k$ is the obtained solution after $k$ iteration of a discrete-algorithm with iterate update \eqref{eq:discrete-update} from initial solution $z_0$, then it holds that $z_k=Z(ks)$ where $Z(ks)$ is the solution at $t=ks$ of the ODE \eqref{eq:infode} with initial solution $Z(0)=z_0$.
Furthermore, the $O(s^r)$-resolution ODE can be viewed as the $r$-th ODE multiscale expansion of \eqref{eq:infode}, and thus its approximation error can be bounded by using multiscale analysis \cite{weinan2011principles}. On the other hand, Theorem \ref{thm:high-resolution} shows that if there exists an ODE that can fully characterize the DTA and $g(z,s)$ is infinitely differentiable in $z$ and $s$, the coefficients of the ODE must be recursively given by \eqref{eq:obtainf}.
\vspace{0.2cm}

The next corollary is an application of Theorem \ref{thm:high-resolution} to the three algorithms -- GDA \eqref{eq:gd}, PPM \eqref{eq:ppm} and EGM \eqref{eq:mp}, which also showcases how to utilize Theorem \ref{thm:high-resolution} to obtain the corresponding order resolution ODEs of a DTA.
\begin{cor}\label{cor:ODEs}
(i) The $O(1)$-resolution ODEs of GDA, PPM and EGM are the same, that is, GF: 
\begin{equation}\label{eq:gf}
    \dZ=-F(Z) \ .
\end{equation}
(ii) The $O(s)$-resolution ODE of GDA is 
\begin{equation}\label{eq:gm_ode}
    \dZ=-F(Z)-\frac{s}{2}\nabla F(Z)F(Z) \ .
\end{equation}
(iii) The $O(s)$-resolution ODEs of PPM and of EGM are the same:
\begin{equation}\label{eq:ppm_ode}
    \dZ=-F(Z)+\frac{s}{2}\nabla F(Z)F(Z) \ .
\end{equation}
\red{
(iv) The $O(s^2)$-resolution ODE of PPM is:
\begin{equation}\label{eq:Os2ppm_ode}
    \dZ=-F(Z)+\frac{s}{2}\nabla F(Z)F(Z) + s^2\pran{-\tfrac{1}{3}(\nabla F(Z))^2 F(Z) - \tfrac{1}{12}\nabla^2 F(Z)(F(Z),F(Z))} \ .
\end{equation}
(v) The $O(s^2)$-resolution ODE of EGM is:
\begin{equation}\label{eq:Os2egm_ode}
    \dZ=-F(Z)+\frac{s}{2}\nabla F(Z)F(Z) + s^2\pran{\tfrac{2}{3}(\nabla F(Z))^2 F(Z) - \tfrac{1}{12}\nabla^2 F(Z)(F(Z),F(Z))} \ .
\end{equation}
}
\end{cor}

\textbf{Proof.}
For GDA with iterate update \eqref{eq:gd}, we have $z^+=z-sF(z)$, thus $g_0(z)=z$, $g_1(z)=-F(z)$ and $g_2(z)=0$ in the Taylor expansion of $g(z,s)$ \eqref{eq:expan_g}. It then follows by the recursive rule \eqref{eq:obtainf} that
\begin{align}
    \begin{split}
        f_0(Z)&=g_1(Z)=-F(Z) \\
        f_1(Z)&=g_2(Z)-\frac{1}{2}h_{2,0}(Z)=0- \frac{1}{2}\nabla f_0(Z)f_0(Z)=-\frac{1}{2}\nabla F(Z)F(Z) \ ,
    \end{split}{}
\end{align}
therefore the $O(1)$-resolution ODE of GDA is \eqref{eq:gf} and the $O(s)$-resolution ODE of GDA is \eqref{eq:gm_ode}.

\red{
For PPM with iterate update \eqref{eq:ppm}, we have $z^+=z-sF(z^+)$, thus by expanding the operator $(I+sF)^{-1}$, we obtain
\begin{align}\label{eq:second-order-Taylor-PPM}
\begin{split}
        z^+ &=g(z,s) = (I+sF)^{-1}(z) \\
        &= z -sF(z) + s^2 \nabla F(z) F(z) + s^3\pran{ -(\nabla F(z))^2 F(z) - \frac{1}{2}\nabla^2 F(z) (F(z), F(z))} + o(s^3) \ ,
\end{split}
\end{align}
whereby $g_0(z)=z$, $g_1(z)=-F(z)$, $g_2(z)=\nabla F(z) F(z)$ and $g_3(z)=-(\nabla F(z))^2 F(z) - \frac{1}{2}\nabla^2 F(z) (F(z), F(z))$ in the Taylor expansion of $g(z,s)$ \eqref{eq:expan_g}, where $\nabla^2 F(z)$ is a tensor and $\nabla^2 F(z) (F(z), F(z))$ refers to tensor product {(For the completeness of the paper, we present the calculation of the expansion \eqref{eq:second-order-Taylor-PPM} in Appendix \ref{sec:expansion_operator})}. It then follows by the logic flow (Figure \ref{fig:logic_flow}) and the recursive rule \eqref{eq:recursive}\eqref{eq:recursive_compute_h} that
\begin{align}\label{eq:cal_PPM_ODE}
    \begin{split}
        f_0(Z)&=h_{1,0}(Z)=-F(Z) \\
        h_{2,0}(Z)&=\nabla h_{1,0}(Z) h_{1,0}(Z) = \nabla F(Z)F(Z)\\
        f_1(Z)&=h_{1,1}(Z)=g_2(Z)-\tfrac{1}{2}h_{2,0}(Z)=\tfrac{1}{2}\nabla F(Z)F(Z) \\
        h_{2,1}(Z) &=\nabla h_{1,0}(Z) h_{1,1}(Z)+\nabla h_{1,1}(Z) h_{1,0}(Z)= -(\nabla F(Z))^2 F(Z)-\tfrac{1}{2}\nabla^2F(Z)(F(Z), F(Z)) \\
        h_{3,0}(Z) &= \nabla h_{2,0}(Z) h_{1,0}(Z) = -(\nabla F(Z))^2 F(Z) - \nabla^2 F(Z)(F(Z), F(Z)) \\
        f_{2}(Z) &= g_3(Z)-\tfrac{1}{2}h_{2,1}(Z)-\tfrac{1}{6}h_{3,0}(Z)= -\tfrac{1}{3}(\nabla F(Z))^2 F(Z) - \tfrac{1}{12}\nabla^2 F(Z)(F(Z),F(Z)) \ ,
    \end{split}{}
\end{align}
therefore the $O(1)$-resolution ODE of PPM is \eqref{eq:gf} and the $O(s)$-resolution ODE of GDA is \eqref{eq:ppm_ode}.

For EGM with iterate update \eqref{eq:mp}, we have
\begin{equation*}
    z^+ = z-sF(z-sF(z)) = z -sF(z) + s^2 \nabla F(z) F(z) - \frac{s^3}{2} \nabla^2 F(z)(F(z), F(z))+ o(s^3) \ ,
\end{equation*}
whereby $g_0(z)=z$, $g_1(z)=-F(z)$, $g_2(z)=\nabla F(z) F(z)$ and $g_3(z)=-\frac{1}{2} \nabla^2 F(z)(F(z), F(z))$ in the Taylor expansion of $g(z,s)$ \eqref{eq:expan_g}. Following the same calculation as \eqref{eq:cal_PPM_ODE}, we have that $f_2(Z)=\tfrac{2}{3}(\nabla F(Z))^2 F(Z) - \tfrac{1}{12}\nabla^2 F(Z)(F(Z),F(Z))$, which finishes the proof.
 \qed
}

\red{
In the end of this section, we highlight that the above $O(s^r)$-resolution ODE framework can be used to analyze generic DTAs with iterate update $g(z,s)$. Some potential applications include but not limited to (i) analyzing other algorithms for minimax problems, such as Alternating Gradient Descent Ascent (AGDA), PDHG~\cite{chambolle2011first} and ADMM~\cite{douglas1956numerical,eckstein1992douglas}, etc; (ii) analyzing continuous optimization methods, such as gradient descent, mirror descent, Newton's method, etc; (iii) finding equilibrium of multi-player finite games when the evolving dynamic is continuous (for example logit response dynamic~\cite{blume1993statistical}). However, this framework does not apply directly to Nesterov's accelerated method for minimizing a strongly-convex function~\cite{nesterovBook}, because $g(z,0)\not = z$ due to the existence of the momentum term in the algorithm, which violates our assumption on the function $g(z,s)$.

}

\vspace{0.2cm}

\subsection{Understanding the Behaviors of DTAs Using Their $O(s^r)$-Resolution ODEs}\label{sec:understanding_behaviors}

In this section, we explain the puzzling behaviors of GDA, PPM, EGM for solving the minimax problems \eqref{eq:POI} via their corresponding ODEs. Informally, we call a certain behavior (such as convergent, divergent, etc) of a DTA an $O(s^r)$-behavior if such behavior can be captured by its $O(s^r)$-resolution ODE. Moreover, if different algorithms correspond to the same $O(s^r)$-resolution ODE, then they should exhibit similar $O(s^r)$-behavior (upto a smaller order difference) from the multi-scale analysis viewpoint \cite{weinan2011principles}. This argument will be formalized later in Section \ref{sec:ODEtoDTA}.

In Scenario (i) when $L(x,y)$ is $\mu$-strongly convex-strongly concave, GF converges linearly to the unique stationary point of \eqref{eq:POI}. This is an $O(1)$-behavior. To see it, we observe that $\|F(Z)\|^2$ is a linear decaying energy function of GF \eqref{eq:GF} \footnote{This type of decaying rate is called ``exponential rate'' in ODE literature. We here use the terminology ``linear rate'' in order to be consistent with the linear convergence in optimization literature.}:
\begin{align*}
    \begin{split}
        \frac{d}{dt} \frac{1}{2}\|F(Z)\|^2 & = F(Z)^T \nabla F(Z) \dZ = -F(Z)^T \nabla F(Z) F(Z)= -F(Z)^T \twomatrix{\nabla_{xx} L(x,y)}{\nabla_{xy} L(x,y)}{-\nabla_{xy}L(x,y)^T}{\nabla_{yy} L(x,y)} F(Z)\\
        & = -F(Z)^T \twomatrix{\nabla_{xx} L(x,y)}{}{}{\nabla_{yy} L(x,y)} F(Z) \le -\mu \|F(Z)\|^2 \ ,
    \end{split}
\end{align*}
thus $\|F(Z(t))\|^2\le \exp(-2\mu t)\|F(Z(0))\|^2$. Notice that the above linear convergence rate of GF is $O(1)$ (since the $2\mu$ term in the linear rate is independent of $s$), and the $O(1)$-resolution ODEs of GDA, PPM and EGM are all GF, which intuitively explains why GDA, PPM and EGM all converge linearly to the solution to \eqref{eq:POI} in Scenario (i) by
following the trajectories as GF. The formal proof of the linear convergence rate of the three discrete-time algorithms in Scenario (i) can be found in~\cite{shi2018understanding,rockafellar1976monotone,tseng1995linear}.

\begin{figure}[!tbp]
\centering
  \begin{subfigure}[b]{0.3\textwidth}
    \includegraphics[width=\textwidth, trim=10cm 0 10cm 0, clip]{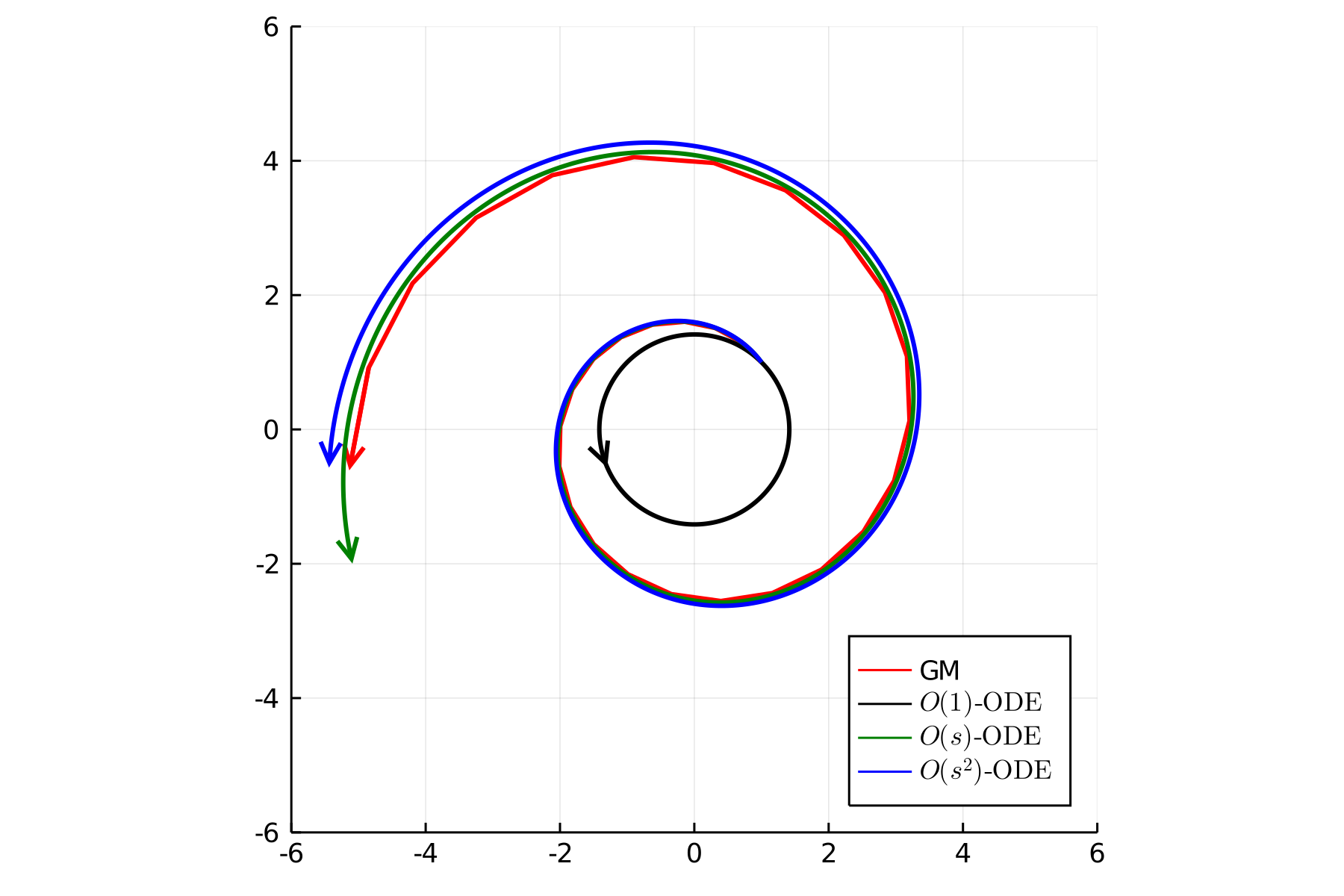}
    \caption{The trajectories of GDA and its corresponding ODEs.}
  \end{subfigure}
  \hspace{0.2cm}
  \begin{subfigure}[b]{0.3\textwidth}
    \includegraphics[width=\textwidth, trim=10cm 0 10cm 0, clip]{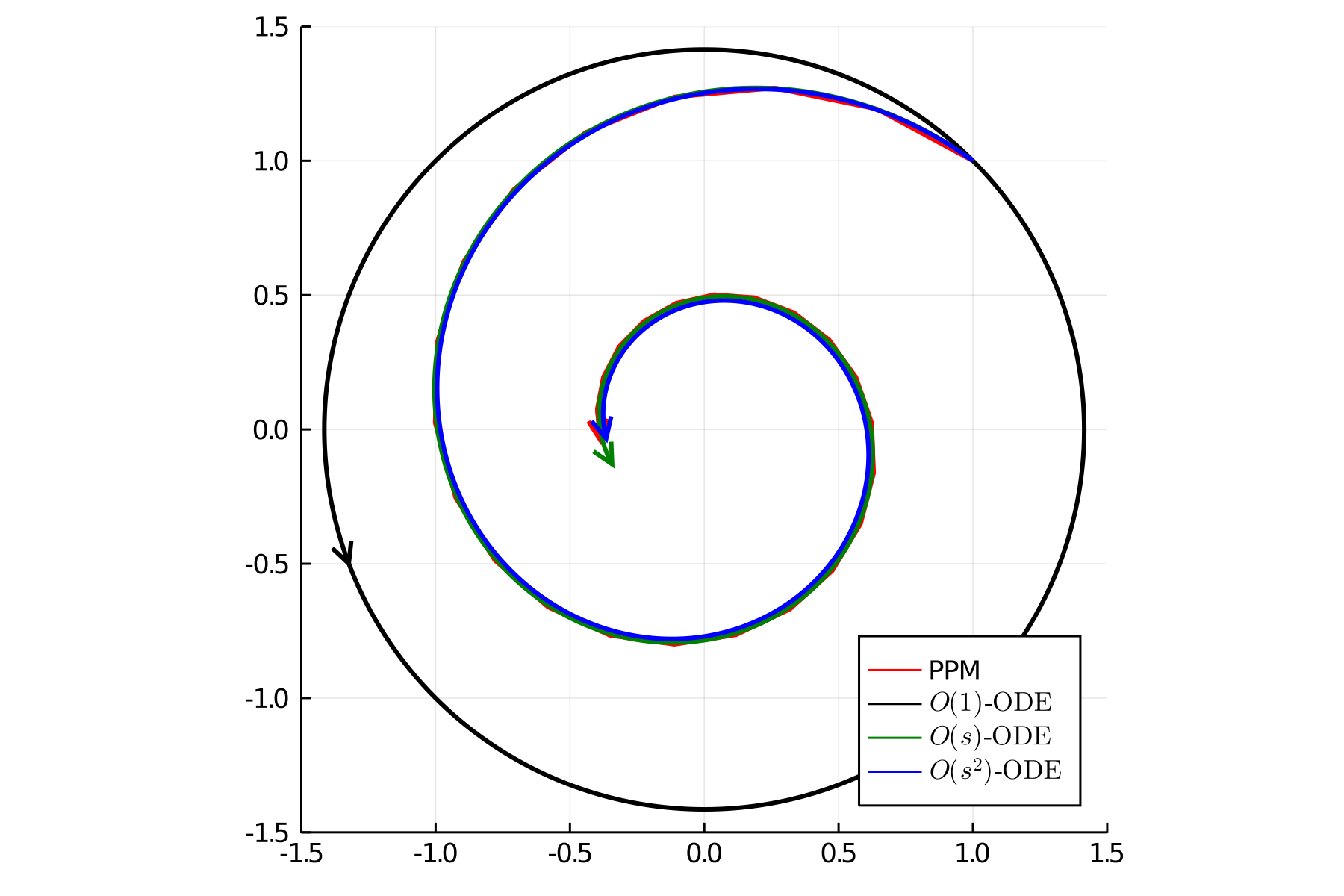}
    \caption{The trajectories of PPM and its corresponding ODEs.}
  \end{subfigure}
  \hspace{0.2cm}
  \begin{subfigure}[b]{0.3\textwidth}
    \includegraphics[width=\textwidth, trim=10cm 0 10cm 0, clip]{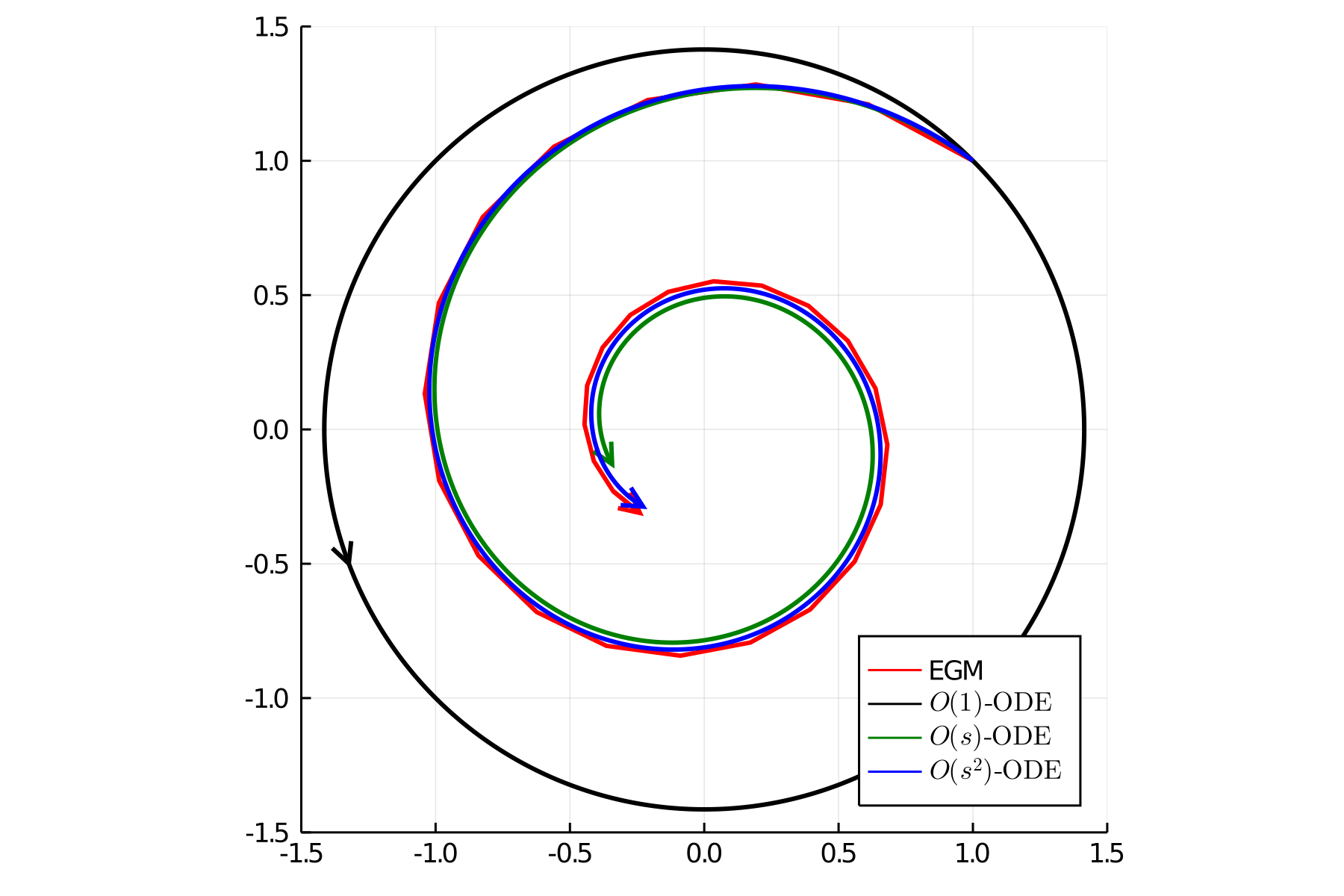}
    \caption{The trajectories of EGM and its corresponding ODEs.}
  \end{subfigure}
  \caption{Illustration of the behaviors of the discrete-time algorithms and their corresponding ODEs. The figure plots the trajectories of different algorithms for solving $\min_x \max_y x y$ with step-size $s=0.3$ and initial solution $(1,1)$.}
  \label{fig:ODE}
\end{figure}

However, the $O(1)$-resolution ODE (i.e. GF \eqref{eq:GF}) does not differentiate between GDA, PPM and EGM, thus it cannot explain the convergent/divergent behaviors of these three algorithms in Scenario (ii). Figure \ref{fig:ODE} plots the trajectories of GDA, PPM and EGM as well as their $O(1)$, $O(s)$ and $O(s^2)$-resolution ODEs in Scenario (ii). As we can see, the higher the order of resolution, the smaller the gap between the trajectory of DTA and the ODE.  Indeed, the convergent/divergent behaviors of GDA, PPM and EGM can be explained with their $O(s)$-resolution ODE as follow (thereby they are $O(s)$-behaviors):

Recall that in Scenario (ii), we consider the bilinear problem
\begin{equation}\label{eq:bilinear}
    \min_{x}\max_y x^T B y \ ,
\end{equation}
thus $F(z)=\twomatrix{}{B}{-B^T}{} z$ and $\nabla F(z)=\twomatrix{}{B}{-B^T}{}$. The $O(s)$-resolution ODE of PPM and EGM \eqref{eq:ppm_ode} becomes
\begin{equation}\label{eq:dzbilinear}
    \dZ=-\twomatrix{}{B}{-B^T}{}Z -\frac{s}{2} \twomatrix{BB^T}{}{}{B^TB} Z =  \twomatrix{-\frac{s}{2}BB^T}{-B}{B^T}{-\frac{s}{2}B^TB} Z.
\end{equation}
Suppose the SVD of $B$ is $B=U^T DV$, where $D$ is an $n$ by $m$ diagonal matrix with $p$ non-zero entries. Then we can rewrite \eqref{eq:dzbilinear} by changing basis $\hZ=\twomatrix{U}{}{}{V}Z$ as
\begin{equation}
    \dot{\hZ}=  \twomatrix{-\frac{s}{2}DD^T}{-D}{D^T}{-\frac{s}{2}D^TD} \hZ.
\end{equation}
Under such basis, there are $p$ independently evolving $2$-d ODE systems, and the $i$-th one is
\begin{align}\label{eq:2dode}
    \begin{split}
        \dhx_i = -\tfrac{s\lambda_i^2}{2} \hx_i - \lambda_i \hy_i\ \  \ , \ \ \dhy_i = -\tfrac{s\lambda_i^2}{2} \hy_i + \lambda_i \hx_i \ ,
    \end{split}
\end{align}
where $\hx_i$ and $\hy_i$ are the variables corresponding to the $i$-th singular-value $\lambda_i$ of matrix $B$. The solution to \eqref{eq:2dode} is given by
\begin{align}\label{eq:solution-os}
    \begin{split}
        \hx_i(t) = c_i e^{-\frac{s}{2}\lambda_i^2 t} \cos(\lambda_i t + \delta_i)\ \  \ , \ \ 
        \hy_i(t) = c_i e^{-\frac{s}{2}\lambda_i^2 t} \sin(\lambda_i t + \delta_i) \ ,
    \end{split}
\end{align}
where $c_i=\sqrt{\hx_i(0)^2 +\hy_i(0)^2}$ and $\delta_i=\arg\tan(\hy_i(0)/\hx_i(0))$ are constants defined by the initial solution. Noticing that the $e^{-\frac{s}{2}\lambda_i^2 t}$ term goes to $0$ linearly as $t\rightarrow\infty$ and the $\cos(\lambda_i t + \delta_i)$ term introduces periodic oscillation in \eqref{eq:solution-os}, which explains the convergent while circling behavior of PPM and EGM in Figure \ref{fig:ODE} (b) (c). Another observation is that when $t$ is large, the 2-d system \eqref{eq:2dode} corresponding to the smallest non-zero singular-value quickly dominates the dynamic, which implies that the oscillation frequency and linear convergence rate is captured by the smallest non-zero singular-value of matrix $B$.

Similarly, the solution of the $O(s)$-resolution ODE  of GDA \eqref{eq:gm_ode} can be characterized after changing basis by
\begin{align*}
    \begin{split}
        \hx_i(t) = c_i e^{\frac{s}{2}\lambda_i^2 t} \cos(\lambda_i t + \delta_i)  \ \  \ , \ \ 
        \hy_i(t) = c_i e^{\frac{s}{2}\lambda_i^2 t} \sin(\lambda_i t + \delta_i) \ .
    \end{split}
\end{align*}
Noticing that the $e^{\frac{s}{2}\lambda_i^2 t}$ term goes to $+\infty$ linearly as $t\rightarrow\infty$. This explains the divergent while circling behavior of GD in Figure \ref{fig:ODE} (a).

\red{
Furthermore, there is a subtle difference between the trajectories of PPM and EGM in the sense that EGM has slightly higher frequency than its $O(s)$-resolution ODE, while PPM has slightly lower frequency than its $O(s)$-resolution ODE. This phenomenon is an $O(s^2)$-behavior, and can be distinguished from their $O(s^2)$-resolution ODEs. Similar to the above arguments, the $O(s^2)$-resolution ODE of PPM results in independent evolving 2-d ODE systems given by
\begin{align*}
    \begin{split}
        \dhx_i = -\tfrac{s\lambda_i^2}{2} \hx_i - \pran{\lambda_i-\tfrac{s^2 \lambda_i^3}{3}} \hy_i \ \  \ , \ \ 
        \dhy_i = -\tfrac{s\lambda_i^2}{2} \hy_i + \pran{\lambda_i-\tfrac{s^2 \lambda_i^3}{3}} \hx_i \ ,
    \end{split}
\end{align*}
whose solutions are:
\begin{align*}
    \begin{split}
        \hx_i(t) = c_i e^{-\frac{s}{2}\lambda_i^2 t} \cos((\lambda_i-\tfrac{s^2}{3}\lambda_i^3) t + \delta_i) \ \  \ , \ \ 
        \hy_i(t) = c_i e^{-\frac{s}{2}\lambda_i^2 t} \sin((\lambda_i-\tfrac{s^2}{3}\lambda_i^3) t + \delta_i) \ .
    \end{split}
\end{align*}
The $-\tfrac{s^2}{3}\lambda_i^3$ term in the frequency explains the lower frequency of PPM compared to its $O(s)$-resolution ODE, as shown in Figure \ref{fig:ODE} (b). In contrast, the corresponding independent evolving 2-d of the $O(s^2)$-resolution ODE of EGM has solutions:
\begin{align*}
    \begin{split}
        \hx_i(t) = c_i e^{-\frac{s}{2}\lambda_i^2 t} \cos((\lambda_i+\tfrac{2s^2}{3}\lambda_i^3) t + \delta_i) \ \  \ , \ \ 
        \hy_i(t) = c_i e^{-\frac{s}{2}\lambda_i^2 t} \sin((\lambda_i+\tfrac{2s^2}{3}\lambda_i^3) t + \delta_i) \ .
    \end{split}
\end{align*}
The $\tfrac{2s^2}{3}\lambda_i^3$ term in the frequency explains the higher frequency of PPM compared to its $O(s)$-resolution ODE, as shown in Figure \ref{fig:ODE} (c).}

\red{
\subsection{Designing New Algorithms Motivated by the $O(s^r)$-Resolution ODEs}\label{sec:design-new-al}

In this section, we present an example to showcase how the $O(s^r)$-resolution ODE framework can help design new optimization algorithms.

From the discussion in the previous section, we know it holds for bilinear minimax problem (i.e., Scenario (ii)) that $\langle-F(Z),Z\rangle=0$, which means the $O(1)$ term in \eqref{eq:gm_ode}\eqref{eq:ppm_ode} is perpendicular to the direction towards the minimax solution, thus it only provides oscillation/circling around the minimax solution. In contrast, the reason PPM/EGM converges while GDA diverges is due to their sign of the $O(s)$ term $\nabla F(Z) F(Z)$, which points directly to the minimax solution. An immediate question is whether we can design a new algorithm that can avoid the oscillation/circling and go directly towards the minimax solution for bilinear minimax problems. A natural idea is to only utilize $O(s)$ term and consider the following ODE: 
\begin{equation}\label{eq:new_ode}
\dZ=\nabla F(Z) F(Z) \ ,    
\end{equation}
whose explicit discretization leads to a new DTA with iterate update
\begin{equation}\label{eq:new_al}
    z^+=z + s \nabla F(z) F(z) \ .
\end{equation}
We call this new algorithm Jacobian method (JM) as it utilizes the Jacobian of $F(z)$. Although JM is a second-order method, it is known that the computational cost of Hessian-gradient product is at the same level of computing the gradient~\cite{pearlmutter1994fast}. Figure \ref{fig:discrete} plots the trajectory of JM. As expected, JM avoids the oscillation and goes toward the minimax solution directly in Figure \ref{fig:discrete} (b). 

Similar to the argument in Section \ref{sec:understanding_behaviors}, we can utilize $O(1)$-resolution ODE \eqref{eq:new_ode} to understand the behaviors of JM \eqref{eq:new_al}. When applying to bilinear problem \eqref{eq:bilinear}, \eqref{eq:new_ode} becomes
\begin{equation}\label{eq:dznewbilinear}
    \dZ=- \twomatrix{BB^T}{}{}{B^TB} Z \ .
\end{equation}
Similar to $O(s)$-resolution ODE of PPM, there are $p$ independent evolving 2-d ODE systems in \eqref{eq:dznewbilinear} after changing basis,
\begin{equation*}
    \dhx_i = -\lambda_i^2 \hx_i \ ,\ \  \dhy_i = -\lambda_i^2 \hy_i \ ,
\end{equation*}
whose solution is given by
\begin{equation*}
        \hx_i(t) = \hx_i(0) e^{-\lambda_i^2 t} \ ,\ \ \        \hy_i(t) = \hy_i(0) e^{-\lambda_i^2 t} \ .
\end{equation*}
Compared with \eqref{eq:solution-os}, we can clearly see that JM avoids oscillations in contrast to the dynamics of PPM and EGM.
}

\section{Analyze the ODE: The $O(s^r)$-Linear-Convergence Conditions}\label{sec:linear-convergence-condition}
In this section, we discuss how to analyze the convergent behavior of the $O(s^r)$-resolution ODE by introducing the $O(s^r)$-linear-convergence condition of a DTA (with respect to an energy function) and presenting examples of such conditions for minimax algorithms.

The typical approach to show that an ODE converges to a fixed point of the dynamic is by identifying an energy function $E(z)$, such that
\begin{itemize}
    \item $E(z(t))$ monotonically decay in $t$;
    \item $E(z)\ge 0$, and $E(z^*)=0$ iff $z^*$ is a fixed point of the dynamic.
\end{itemize}
The convergence of the ODE then can be characterized by the decay rate of the energy function.

We say a condition an $O(s^r)$-linear-convergence condition of a DTA with respect to an energy function $E(z)$ if such condition can guarantee the $O(s^r)$-resolution ODE of a DTA has linear convergence in $E(z)$. More formally,

\begin{mydef}\label{def:linear-conv-condition}
Consider the $O(s^r)$-resolution ODE of a DTA: $\dZ=f^{(r)}(Z,s)$. Suppose there exists a condition which can guarantee that there exists $\rho(s)>0$ such that it holds for any $Z$
\begin{equation}\label{eq:linear-decay}
    \frac{d}{dt} E(Z) = \langle \nabla E(Z), f^{(r)}(Z,s) \rangle \le -\rho(s) E(Z) \ ,
\end{equation}
then we call this condition an $O(s^r)$-linear-convergence condition of the DTA.
\end{mydef}

Inequality \eqref{eq:linear-decay} guarantees that the energy $E(Z)$ decays linearly to $0$ because it holds from \eqref{eq:linear-decay} that $E(Z(t))\le e^{-\rho(s)t}E(Z(0))$. Of course, how to select a good energy function for a specific DTA can be a non-trivial task, and we defer the discussions on this topic in Section \ref{sec:ODEtoDTA}. Here we focus on the inverse problem, that is, given an energy function, we study under what conditions the $O(s^r)$-resolution ODE does have linear convergence.


To further illustrate the idea of the $O(s^r)$-linear-convergence condition, we here present the corresponding conditions of PM, EGM, PPM and JM with energy function:
\begin{equation}\label{eq:continuous-energy}
    E(z)=\frac{1}{2}\|F(z)\|^2\ . \ \ \ \ 
\end{equation}
First, we introduce some new notations that will be used in this section: Denote $A(z)=\nabla_{xx}L(x,y)$, $B(z)=\nabla_{xy}L(x,y)$, $C(z)=-\nabla_{yy}L(x,y)$, then $\nabla F(z)=\left[\begin{smallmatrix} A(z) & B(z) \\ -B(z)^T & C(z) \end{smallmatrix}\right]$. We also use $A, B, C$ to represent $A(z), B(z), C(z)$ if they do not cause any misunderstandings. 
Then
\begin{prop}\label{prop:Oscondition}
(i) An $O(1)$-linear-convergence condition of PM, EGM and PPM is strong convexity-concavity of $L(x,y)$, i.e., there exists $\rho>0$ such that
\begin{equation}\label{eq:O1condition}
    F(Z)^T \left[\begin{matrix} A & 0 \\ 0 & C\end{matrix}\right] F(Z)\ge \frac{\rho}{2}\|F(Z)\|^2\ ,  \text{ for any }Z\ . 
\end{equation}
(ii) An $O(s)$-linear-convergence condition of EGM and PPM is 
\begin{equation}\label{eq:OsconditionPPM}
    F(Z)^T \left[\begin{matrix} A- \frac{s}{2} A^2+ \frac{s}{2} BB^T & 0 \\ 0 & C- \frac{s}{2} C^2+ \frac{s}{2} B^T B \end{matrix}\right] F(Z)\ge \frac{1}{2}\rho(s)\|F(Z)\|^2\ ,  \text{ for any }Z\ .
\end{equation}
(iii) An $O(s)$-linear-convergence condition of GDA is 
\begin{equation}\label{eq:OsconditionGM}
    F(Z)^T \left[\begin{matrix} A+ \frac{s}{2} A^2- \frac{s}{2} BB^T & 0 \\ 0 & C+ \frac{s}{2} C^2- \frac{s}{2} B^T B \end{matrix}\right] F(Z)\ge \frac{1}{2}\rho(s)\|F(Z)\|^2\ ,  \text{ for any }Z\ . 
\end{equation}
\red{
(iv) An $O(1)$-linear-convergence condition of JM is
\begin{equation}\label{eq:O1conditionJM}
    F(Z)^T \left[\begin{matrix} BB^T -A^2 & 0 \\ 0 & B^T B -C^2 \end{matrix}\right] F(Z) \ge \frac{\rho}{2}\|F(Z)\|^2\ .
\end{equation}}
\end{prop}

\textbf{Proof.} 
(i) Substituting the $O(1)$-resolution ODE of GDA, EGM and PPM, namely $\dZ=-F(Z)$, into \eqref{eq:linear-decay}, we obtain 
\begin{align*} 
    \begin{split}
        \frac{d}{dt} \frac{1}{2}\|F(Z)\|^2 & = -F(Z)^T \nabla F(Z) F(Z)= -F(Z)^T \twomatrix{A}{B}{-B^T}{C} F(Z)\\
        & = -F(Z)^T \twomatrix{A}{}{}{C} F(Z) \le -\frac{\rho}{2} \|F(Z)\|^2 \ ,
    \end{split}
\end{align*}
which shows that \eqref{eq:O1condition} is an $O(1)$-linear convergence condition of GDA, EGM and PPM.

(ii) Substituting the $O(s)$-resolution ODE of EGM and PPM, namely $\dZ=-F(Z)+\frac{s}{2}\nabla F(Z)F(Z)$, into \eqref{eq:linear-decay}, we obtain,
\begin{align}\label{eq:LypOs}
    \begin{split}
        \frac{d}{dt} \frac{1}{2}\|F(Z)\|^2 & = F(Z)^T \nabla F(Z) \dZ\\
        & = -F(Z)^T \nabla F(Z) F(Z) + \frac{s}{2} F(Z)^T (\nabla F(Z))^2 F(Z) \\
        & = -F(Z)^T \left[\begin{matrix} A- \frac{s}{2} A^2+ \frac{s}{2} BB^T & 0 \\ 0 & C- \frac{s}{2} C^2+ \frac{s}{2} B^T B \end{matrix}\right] F(Z) \\
        &\le  -\frac{\rho(s)}{2} \|F(Z)\|^2 \ , 
    \end{split}
\end{align}
which shows that \eqref{eq:OsconditionPPM} is an $O(s)$-linear convergence condition of EGM and PPM.

(iii) The proof is the same as (ii) by replacing the sign of the corresponding terms to $\frac{s}{2} F(Z)^T (\nabla F(Z))^2 F(Z)$ in \eqref{eq:LypOs}. 

\red{
(iv) Notice that \eqref{eq:new_ode} is the $O(1)$-resolution ODE of JM. Substituting \eqref{eq:new_ode} into \eqref{eq:linear-decay}, we obtain
\begin{align*}\label{eq:LypOsJM}
    \begin{split}
        \frac{d}{dt} \frac{1}{2}\|F(Z)\|^2 & = F(Z)^T \nabla F(Z) \dZ = F(Z)^T (\nabla F(Z))^2 F(Z) \\
        & = -F(Z)^T \left[\begin{matrix} BB^T -A^2 & 0 \\ 0 & B^T B -C^2 \end{matrix}\right] F(Z)
        \le  -\frac{\rho}{2} \|F(Z)\|^2 \ , 
    \end{split}
\end{align*}
which shows that \eqref{eq:O1conditionJM} is an $O(1)$-linear convergence condition of JM.}
\qed

In the following, we comment on the corresponding linear-convergence conditions of the four algorithms as stated above.

\textbf{($O(s)$-condition of PPM/EGM)} When the step-size $s\le \frac{1}{\lambda}$, a stronger $O(s)$-linear-convergence condition of EGM and PPM for convex-concave problem is
\begin{equation}\label{eq:O(s)-strong}
    F(Z)^T \left[\begin{matrix}  A+ s BB^T & 0\\0  & C+ s B^T B \end{matrix}\right] F(Z)\ge \rho(s)\|F(Z)\|^2\ ,  \text{ for any }Z\ ,
\end{equation}
by noticing $A- \frac{s}{2} A^2+ \frac{s}{2} BB^T \ge \frac{1}{2}(A+ s BB^T)$.
This stronger condition clearly shows that the interaction terms help the linear convergence of EGM and PPM, and in contrast, the interaction terms hurt the linear convergence of GDA, which provides another explanation to the convergent/divergent behaviors of different algorithms in Figure \ref{fig:discrete} (b) when the objective is bilinear. This is consistent with the argument in \cite{liang2019interaction}. Moreover, in this case, the linear rate $\rho(s)$ usually is a linear function in $s$ with nonnegative slope and intercept. Finally we comment that the $O(s)$-resolution ODE of PPM and EGM does not require convexity-concavity of $L(x,y)$ (as long as it has a sufficient interaction term), which is consistent with the recent results on the landscape of PPM for solving nonconvex-nonconcave minimax problems \cite{grimmer2020landscape}.
\red{

\textbf{($O(1)$-condition of GDA and the step-size upper bounds)} It is well-known that under their $O(1)$-linear-convergence condition (i.e., when $L(x,y)$ is strongly-convex-strongly-concave), GDA needs to take smaller step-size ($s\le O(\frac{\mu}{\gamma^2})$) than that for convex optimization ($s\le O(\frac{1}{\gamma})$) in order to obtain linear convergence ~\cite{du2017stochastic}\cite{du2019linear}\cite{grimmer2020landscape}. The reason for the smaller step-size can be clearly seen from their $O(s)$-linear-convergence conditions. Informally speaking, the $O(s)$-linear-convergence condition of GDA \eqref{eq:OsconditionGM} requires $A\succeq O(sBB^T)$ and $C\succeq O(sB^T B)$. A sufficient condition to guarantee that is $s\le O(\frac{\mu}{\gamma^2})$ because if it holds, we have
$A\succeq \mu I \succeq O(s \gamma^2 I) \succeq O(s BB^T)$ (same argument applies to $C$). 


\textbf{($O(1)$-condition of JM)} The $O(1)$-linear-convergence condition of JM \eqref{eq:O1conditionJM} holds when $L(x,y)$ has sufficient interaction term (i.e., $BB^T\succeq A^2$ and $B^TB\succeq C^2$), and such condition may hold for nonconvex-nonconcave minimax problems. 
}
\vspace{0.3cm}

Now we focus on the $O(s)$-linear-convergence condition of PPM and EGM \eqref{eq:OsconditionPPM} in order to study the linear convergence of these two algorithms beyond the two classic scenarios when $L(x,y)$ is either strongly convex-strongly concave or bilinear. Indeed, the condition \eqref{eq:OsconditionPPM}, as well as its weaker version \eqref{eq:O(s)-strong}, is a general condition that is satisfied by many objective $L(x,y)$, and we herein present some examples:

Section \ref{sec:high-resolution} utilizes the corresponding ODE systems of GDA, PPM and EGM to explain their behaviors for solving minimax problem \eqref{eq:POI} in the two classic scenarios when $L(x,y)$ is either strongly convex-strongly concave or bilinear. In this section, we study general minimax function $L(x,y)$ beyond these two classic scenarios. Indeed, the $O(s)$-resolution ODE of PPM and EGM \eqref{eq:ppm_ode} inspire us to introduce the $O(s)$-linear-convergence condition of the two algorithms, and we will show that this condition is well satisfied in general by examples.

\begin{exam} \label{ex:1} Suppose $L(x,y)$ is $\mu$-strongly convex-strongly concave, then it is straight-forward to see that $\rho(s)\ge \mu$. This is Scenario (i) in previous sections.
\end{exam}
\begin{exam}\label{ex:2}
Suppose $L(x,y)=x^T B y$ is a bilinear function, then $\rho(s)=s\lambda_{\min}^+ (BB^T)$ by noticing $\FF=\emph{Range}(B)\times\emph{Range}(B^T)$. This is Scenario (ii) in the previous sections.
\end{exam}
\begin{exam}\label{ex:3}
Suppose $L(x,y)=f(x)+x^T B y - g(y)$ where $f(x)$ is $\mu$-strongly convex in $x$, $g(y)$ is concave in $y$ and $B$ has full column rank, then it holds that $\rho(s)\ge \min\{\mu, s\lambda_{\min} (BB^T)\}$. Actually, a recent work \cite{du2019linear} shows that GDA has a linear convergence rate in this case, and our results in Section \ref{sec:convergence} show that PPM and EGM also exhibit linear convergence in this case.
\end{exam}{}
\begin{exam}\label{ex:4}
Suppose $L(x,y)$ satisfies for any $(x,y)\in\RRmn$ that $\nabla_{xy} L(x,y) $ is square (thus $m=n$) and full rank, and there exists a positive $\mu>0$ such that $$\lambda_{\min} (\nabla_{xy} L(x,y)^T\nabla_{xy} L(x,y)) \ge \mu > 0, \ \ \forall (x,y)\ .$$ Then $\rho(s)\ge s \mu$. A more specific example can be $L(x,y)=f(x)+x^T B y - g(y)$ with square and full-rank matrix $B$.
\end{exam}
\begin{exam}\label{ex:5}
Suppose $L(x,y)= f(C_1 x) + x^T B y - g(C_2 y)$ where $f(\cdot)$ and $g(\cdot)$ are both strongly convex. Then we can show that $L(x,y)$ satisfies the $O(s)$-linear-convergence condition \eqref{eq:O(s)-strong} with $\rho(s)>0$. We leave the definition of $\rho(s)$ and the proof of this example in Appendix \ref{sec:app-example}. 
\end{exam}
\begin{exam}\label{ex:6}
Suppose $L(x,y)$ is nonconvex-nonconcave but has a sufficient interaction term such that \eqref{eq:OsconditionPPM} is satisfied. 
\end{exam}

\begin{rem}
Example \ref{ex:3}, \ref{ex:4}, \ref{ex:5}, \ref{ex:6} and the results in Section \ref{sec:ODEtoDTA} show that PPM and EGM have linear convergence for solving \eqref{eq:POI} beyond the two standard scenarios.
\end{rem}

\section{From ODEs back to DTAs: Proper Energy Functions}\label{sec:ODEtoDTA}
Section \ref{sec:high-resolution} presents how to obtain a suitable ODE from a DTA. Section \ref{sec:linear-convergence-condition} presents how to analyze the corresponding ODEs by introducing the $O(s^r)$-linear-convergence condition of a DTA with respect to an energy function. In this section, we close the loop by building up the connections between the convergence of the DTA and its $O(s^r)$-resolution ODE. Informally speaking, we show that with a \emph{proper} choice of the energy function, the linear convergence of its $O(s^r)$-resolution ODE can automatically guarantee that the DTA converges at the same linear convergence rate.


To study the connection between a DTA and its $O(s^r)$-resolution ODE, we begin with discussing the relationship between their fixed points, defined as:
\begin{mydef}
Consider a DTA with iterate update $z^+=g(z,s)$ and its $O(s^r)$-resolution ODE $\dZ=f^{(r)}(Z,s)$ \eqref{eq:rthode}. 

1. We say $z^*$ is a fixed point of the DTA if there exists $s^*>0$ such that $g(z^*,s)=z^*$ for any step-size $s\in (0, s^*]$. 

2. We say $z^*$ is a fixed point of the $O(s^r)$-resolution ODE if there exists $s^*>0$ such that $f^{(r)}(z^*, s)=0$ for step-size $s\in (0, s^*]$.
\end{mydef}
The next proposition connects the fixed points of the DTA and its $O(s^r)$-resolution ODEs.
\begin{prop}\label{prop:fix-point}
Consider a DTA with iterate update $z^+=g(z,s)$ and its $O(s^r)$-resolution ODE $\dZ=f^{(r)}(Z,s)$ \eqref{eq:rthode}.
\begin{enumerate} 
    \item Suppose $z^*$ is a fixed point of the DTA, then $z^*$ is also a fixed point of the $O(s^r)$-resolution ODE for any degree $r$.
    \item Suppose $z^*$ is a fixed point of the $O(s^r)$-resolution ODE of a DTA. Then $g_j(z^*)=0$ for $j=0,...,r+1$, where $g_j$ is the $j$-th coefficient function in the Taylor expansion of $g(z,s)$ (see \eqref{eq:expan_g}).
\end{enumerate}
\end{prop}
\textbf{Proof.} 
1. We prove the claim by contradiction. Consider the $O(s^r)$-resolution ODE \eqref{eq:rthode} to a DTA $z^+=g(z,s)$. If the claim does not hold, then there exists $j\le r$, such that $f_j(z^*)\not=0$, and without loss of generality, let $j$ be the smallest term that $f_j(z^*)\not=0$. Then we know $z^*$ is not a fixed point of the $O(s^j)$-resolution ODE, because $f^{(j)}(z^*,s)$ is a $j$-th degree polynomial in $s$ with at most $j$ different roots.
Thus, it follows from the ODE $\dZ=f^{(j)}(z^*,s)$ that
\begin{equation*}
    \|Z(s)-z^+\|=\|Z(s)-z^*\|= \|Z(s)-Z(0)\|=\Omega(\|s^{j+1} f_j(z^*)\|) \ , \text{ when } s\rightarrow 0 \ \footnote{Recall the $\Omega$ notation means that there exists a constant $c>0$ such that $\|Z(s)-z^+\|\ge c\|s^{j+1} f_j(z^*)\|$ as $s\rightarrow 0$.},
\end{equation*}
where $Z(s)$ is the solution obtained at $t=s$ following the $O(s^j)$-resolution ODE with initial solution $Z(0)=z^*$. This contradicts with the definition of the $O(s^j)$-resolution ODE \eqref{eq:small_term}.

2. Notice that it follows from the definition of the fixed point of the $O(s^r)$-resolution ODE that $f_0(z^*)=f_1(z^*)=...=f_r(z^*)=0$, thus $Z(t)=z^*$ for any $t\ge 0$ following ODE \eqref{eq:rthode} with initial solution $Z(0)=z^*$. The claim follows directly by noticing $h_{j,i}(z^*)=0$ for $i=0,...r-1$ from \eqref{eq:djdt}, thus $g_j(z^*)$=0 from \eqref{eq:obtainf}. 
\qed

Indeed, for many optimization algorithms, in particular first-order methods, $g_1(z^*)=0$ implies $z^*$ is a fixed point of the DTA. This is because, for first-order methods, such as PPM, EGM, GDA discussed in the paper, $g_1(z)$ is usually the gradient of the objective function (upto a scalar). In such a case, Proposition \ref{prop:fix-point} shows the equivalence of the fixed points of these DTAs and its corresponding $O(s^r)$-resolution ODE (for any degree $r$).

Although the fixed points of the DTA and the ODEs are in many cases the same, the linear convergence of the $O(s^r)$-resolution ODE itself, unfortunately, is not enough to guarantee the linear convergence of the DTA. To bridge such a gap, we introduce the \emph{properness} of an energy function that is used in the linear convergence argument for the ODE:
\begin{mydef}\label{def:proper_energy}
We say an energy function $E(z)=\frac{1}{2} e(z)^2$ is proper to study the $O(s^r)$-resolution ODE of a DTA if there exists $c>0$ such that it holds for any $\delta\ge 0$ and $z\in\{e(z) \le \delta\}$ that
\begin{equation}\label{eq:gradient-based-ODE}
    \|Z(s)-z^+\|\le c s^{r+2}  e(z) \ ,
\end{equation}
where $z^+=g(z,s)$ is the output of the DTA from $z$, and $Z(s)$ is the solution obtained at $t=s$ following the $O(s^r)$-resolution ODE \eqref{eq:rthode} with initial solution $Z(0)=z$.
\end{mydef}

Recall that the $O(s^r)$-resolution ODE guarantees that $\|Z(s)-z^+\|\le O(s^{r+2})$ (See Remark \ref{rem:high-order-difference}). Proper energy function \eqref{eq:gradient-based-ODE} further imposes an upper bound on the one-iteration gap $\|Z(s)-z^+\|$ in terms of $z$. A proper energy function always exists, because we can always set $e(z)=\frac{\|Z(s)-z^+\|}{s^{r+2}}$, where we utilize the fact that $\|Z(s)-z^+\|=O(s^{r+2})$ so that $e(z)$ does not blow up as $s\rightarrow 0$, and the fact that $e(z^*)=0$ by noticing $Z(s)=z^+=z^*$ with initial solution $z=z^*$. 

Meanwhile, in order to obtain a more meaningful $O(s^r)$-linear-convergence condition as stated in Section \ref{sec:linear-convergence-condition}, we prefer a simple form of $e(z)$. Some typical examples of $e(z)$ include:
\begin{itemize}
    \item Norm of gradient, i.e., $\|F(z)\|$;
    \item Distance from the current iterate to optimal solutions, i.e., $\|z-z^*\|$;
    \item Square root of the optimality gap for convex optimization;
    \item Linear combination of the above.
\end{itemize}

\vspace{0.2cm}
    
    


%



Let $S^0=\{z|E(z)\le E(z^0)\}$ be the level set of $E$. The next theorem presents our main result that bridges the convergence of a DTA and its $O(s^r)$-resolution through a proper energy function:

\begin{thm}\label{thm:gradient-based-ODE-bound}
Consider a DTA with iterate update $z^+=g(z,s)$ and its $O(s^r)$-resolution ODE $\dZ=f^{(r)}(Z,s)$. Suppose 

(i) the $O(s^r)$-resolution ODE converges to an optimal solution with respect to a proper energy function, namely, \eqref{eq:linear-decay} holds with a proper energy function $E$; 

(ii) there exists a constant $\gamma$ such that $\|\nabla e(z)\|\le \gamma$ for any $z\in \text{Conv}\pran{S^0\cup\{g(z,s)|z\in S^0\}}$, where $\text{Conv}(\cdot, \cdot)$ denotes the convex hull of two sets; 

(iii) the step-size $s$ satisfies
\begin{equation}\label{eq:step-size-condition}
    \gamma cs^{r+2}\le \min\pran{1, \frac{s\rho(s)}{16}} \ ,
\end{equation}
where $c$ is from the properness of the energy function \eqref{eq:gradient-based-ODE} when choosing $\delta=e(z^0)$. 

Then it holds for any $k\ge 0$ that
\begin{equation*}
    E(z^k)\le \pran{1-\frac{s\rho(s)}{4}}^k E(z^0) \ .
\end{equation*}
\end{thm}


\textbf{Proof.} It follows from Taylor expansion of $E(z)$ that
\begin{align}
    \begin{split}
        E(z^+) & = E(Z(s))+\int_{0}^1 \nabla E(Z(s)+t(z^+-Z(s))) (z^+-Z(s)) dt \\
        &\le E(Z(s)) + \gamma \|z^+-Z(s)\| \int_{0}^1 e(Z(s)+t(z^+-Z(s)))  dt  \\
        & \le E(Z(s)) + \gamma \|z^+-Z(s)\| \int_{0}^1 e(Z(s))+ \gamma \|t(z^+-Z(s))\|  dt \\
         & \le e^{-s\rho(s)} E(z) + \gamma cs^{r+2} e(z) \pran{e(z)+\frac{\gamma}{2} cs^{r+2}e(z)} \\
         & \le \pran{1-\frac{s\rho(s)}{2}} E(z) + 4\gamma cs^{r+2} E(z) \\
         & \le \pran{1-\frac{s\rho(s)}{4}} E(z) \ ,
    \end{split}
\end{align}
where the first inequality utilizes $\|\nabla E(z)\|=\|\nabla e(z) e(z)\|\le \gamma e(z)$, the second inequality utilizes (ii), the third inequality is due to  \eqref{eq:gradient-based-ODE} and \eqref{eq:linear-decay}, and the last two inequality utilizes \eqref{eq:step-size-condition}. This finishes the proof by telescoping. \qed

\reddd{
\begin{rem}
We here examine the three conditions stated in the theorem. (i) requires the energy function $E$ is proper with respect to the $O(s^r)$-resolution ODE. We will present a simple approach to check whether an energy function is proper later in Theorem \ref{thm:gradient-based-ODE}. (ii) requires $e(z)$ to be Lipschitz continuous in set $\text{Conv}\pran{S^0, \{g(z,s)|z\in S^0\}}$. In many cases, the level set $S^0$ is close and bounded, so as $\text{Conv}\pran{S^0, \{g(z,s)|z\in S^0\}}$, thus (ii) is naturally satisfied. In our examples, $e(z)$ is often chosen as distance to the optimal solutions $\|z-z^*\|$ or norm of gradient $\|F(z)\|$, where (ii) is satisfied globally for the former with $\gamma=1$, and for latter when the gradient $F(z)$ is Lipschitz continuous. For (iii), recall that $\rho(s)$ (defined in \eqref{eq:linear-decay}) is usually an $r$-th order polynomial on $s$ with non-negative coefficients due to the construction of the $O(s^r)$-linear-convergence condition (see Proposition \ref{prop:Oscondition} for examples). In such a case, \eqref{eq:step-size-condition} holds with reasonably a small step-size $s$. We present examples of DTAs that satisfy such conditions in Corollary \ref{cor:ODE-to-DTA}. Furthermore, the maximal step-size that guarantees linear-convergence depends on the value of $c$, which we will revisit later in Remark \ref{rem:step-size}.
\end{rem}
}

\medskip

Notice that to verify whether an energy function is proper from definition \eqref{eq:gradient-based-ODE} requires to solve the $O(s^r)$-resolution ODE, which can be highly nontrivial. To avoid this, Theorem \ref{thm:gradient-based-ODE} presents easy-to-check sufficient conditions for proper energy functions. Roughly speaking, if $\|f_j(z)\|$ (or $\|g_j(z)\|$) is upper bounded by $e(z)$, and its high order derivatives are bounded for $z\in S^0$, then the energy function is proper.
\begin{thm}\label{thm:gradient-based-ODE}
Consider the $O(s^r)$-resolution ODE \eqref{eq:rthode} of a DTA with Taylor expansion \eqref{eq:expan_g} and step-size $s<1$. Suppose for any $\delta$ and $z\in \{z|e(z)\le \delta\}$, there exists a constant $a>0$ such that it holds
\begin{equation}\label{eq:gradient-based}
    \|z^+-z\|\le a s e(z) \ ,
\end{equation}
and $\gamma=\max_{z\in S}\|\nabla e(z)\|<\infty$, where $$S:=\text{Conv}(\{g(z,t)|0\le t\le s, z\in \{z|e(z)\le e(z^0)\}\})\ .$$
Suppose either of the following two conditions hold:

(i) \textbf{(conditions on $f_j(z)$)} $f_j(z)$ is $(r+1)$-th order differentiable, and it holds for any $z\in S$ that $$\|f_j(z)\| \le O(e(z)) \text{ and } \|\nabla^k f_j(z)\| \le O(1) \text{ for } j=0,...,r+1 \text{ and } k=1,...,r+1 \ ;$$

(ii) \textbf{(conditions on $g_j(z)$)} $g_j(z)$ is $(2r+3-j)$-th order differentiable over $z$, and it holds for any $z\in S$ that
\begin{equation}\label{eq:condition-on-g}
    \|g_j(z)\| \le O(e(z)) \text{ and } \|\nabla^k g_j(z)\| \le O(1) \text{ for } j=1,...,r+2 \text{ and } k=1,...,2r+3-j \ .
\end{equation}


Then the energy function $E(z)=\frac{1}{2} e(z)^2$ is proper to study the $O(s^r)$-resolution ODE.
\end{thm}

\reddd{
\begin{rem}
We here comment on the implication of Theorem \ref{thm:gradient-based-ODE}. In order to make sure the gap between one iteration of the DTA and the ODE is upper-bounded by $e(z)$ (namely \eqref{eq:gradient-based-ODE} holds), it is not surprising that we require the movement of one iteration of the DTA is upper-bounded by $e(z)$ (namely \eqref{eq:gradient-based} holds). Moreover, \eqref{eq:gradient-based} is easy to check since it is a condition on the DTA (not the ODE). Meanwhile, notice $S$ is usually a closed and bounded set, in particular when the optimal solution set is bounded, in which case $\|\nabla e(z)\|$  and $\|\nabla^k f_j(z)\|$ (or $\|\nabla^k g_j(z)\|$) is upper bounded for $z\in S$. The most important conditions required in Theorem \ref{thm:gradient-based-ODE} is $\|f_j(z)\| \le O(e(z))$ (or $\|g_j(z)\| \le O(e(z))$), and the critical region is when $z$ is close to an optimal solution thus $e(z)$ is small. In other words, in order to make sure $E(z)=\frac{1}{2}e(z)^2$ is a proper energy function, we essentially require $e(z)$ to be able to upper bound $\|f_j(z)\|$ (or $\|g_j(z)\|$) as $z$ goes to a fixed point $z^*$.
\end{rem}

Now we have all pieces needed in the $O(s^r)$-resolution ODE framework.
As applications to Theorem \ref{thm:gradient-based-ODE-bound} and Theorem \ref{thm:gradient-based-ODE}, the following corollary shows that GDA, PPM, EGM and JM converge linearly to a minimax solution under the corresponding linear-convergence-condition when we choose the energy function $E(z)=\frac{1}{2}\|F(z)\|^2$: 
\begin{cor}\label{cor:ODE-to-DTA}
Denote $S:=\text{Conv}(\{g(z,t)|0\le t\le s, z\in \{z|\|F(z)\|\le \|F(z^0)\|\}\})$.

(i) Suppose $L(x,y)$ is third-order differentiable and $
\|\nabla^j F(z)\|$ is bounded for $j=1,2$ and $z\in S$. Suppose the $O(1)$-linear-convergence condition of GDA, PPM and EGM \eqref{eq:O1condition} holds with $\rho>0$. Then there exists $s^*$ such that for any $s\le s^*$, GDA, PPM and EGM converge linearly to a stationary point of $L(x,y)$.

(ii) Suppose $L(x,y)$ is fifth-order differentiable, and $\|\nabla^j F(z)\|$ is bounded for $j=1,\ldots, 4$ and $z\in S$. Suppose the $O(s)$-linear-convergence condition of PPM and EGM \eqref{eq:OsconditionPPM} holds with $\rho(s)\ge ds$ for $d>0$ and small $s$. Then there exists $s^*$ such that for any $s\le s^*$, PPM and EGM converge linearly to a stationary point of $L(x,y)$.

(iii) Suppose $L(x,y)$ is fifth-order differentiable, and $\|\nabla^j F(z)\|$ is bounded for $j=1,\ldots, 4$ and $z\in S$. Suppose the $O(s)$-linear-convergence condition of GDA \eqref{eq:OsconditionGM} holds with $\rho(s)\ge ds$ for $d>0$ and small $s$. Then there exists $s^*$ such that for any $s\le s^*$, GDA converges linearly to a stationary point of $L(x,y)$.

(iv) Suppose $L(x,y)$ is fourth-order differentiable, and $\|\nabla^j F(z)\|$ is bounded for $j=1,2,3$ and $z\in S$. Suppose the $O(1)$-linear-convergence condition of JM \eqref{eq:O1conditionJM} holds. Then there exists $s^*$ such that for any $s\le s^*$, JM converges linearly to a stationary point of $L(x,y)$.
\end{cor}

\textbf{Proof.} Here we just show (ii) for PPM, and the other claims follow with a similar argument. Recall that $g_0, g_1, g_2, g_3$ for PPM is defined in \eqref{eq:second-order-Taylor-PPM}. Then, it is easy to check that $g_j$ is $(4-j)$-th order differentiable and \eqref{eq:condition-on-g} holds by utilizing the continuity conditions stated in (ii). Furthermore, it holds that $$\|z^+-z\|=s\|F(z^+)\|\le se(z)+s\gamma \|z^+-z\|\ ,$$ thereby $\|z^+-z\|\le \frac{1}{1-s\gamma} s e(z)$. Thus, \eqref{eq:gradient-based} holds with $a=\frac{1}{1-s\gamma}$. It then follows from Theorem \ref{thm:gradient-based-ODE} that the energy function $E(z)=\frac{1}{2}\|F(z)\|^2$ is proper to study the $O(s)$-resolution ODE of PPM (i.e. \eqref{eq:ppm_ode}), thus (i) in Theorem 2 holds. Furthermore, notice $e(z)=\|F(z)\|$, thus $\|\nabla e(z)\|\le \|\nabla F(z)\|$ is bounded in $S
\supseteq \text{Conv}\pran{S^0\cup\{g(z,s)|z\in S^0\}}$, thus (ii) in Theorem 2 holds. Moreover, $\rho(s)$ is a linear function in $s$ with non-negative coefficients (see Proposition \ref{prop:Oscondition} and the discussions afterwards), thus there exists $s^*$ such that \eqref{eq:step-size-condition} holds for any $s\le s^*$. Therefore, it follows from Theorem \ref{thm:gradient-based-ODE-bound} that $E(z^k)$ decays to $0$ linearly, which showcases the linear convergence of PPM. \qed

\begin{rem}
When the level set $S$ is bounded, Corollary \ref{cor:ODE-to-DTA} shows that as long as $F(z)$ is sufficiently differentiable, then the $O(s)$-linear-convergence condition is sufficient to guarantee the linear convergence of GDA, PPM, EGM and JM.
\end{rem}
}
\medskip



The next proposition will be used in the proof of Theorem \ref{thm:gradient-based-ODE}. 

\begin{prop}\label{prop:bound_on_h}
Under either condition stated in Theorem \ref{thm:gradient-based-ODE}, it holds for $j=1,...r+2$ and $i=0,1,...,r(r+2)$ that
$$
h_{j,i}(z) \le O(e(z)) \ ,
$$
where $h_{j,i}(z)$ is defined in \eqref{eq:djdt}. Furthermore, it holds for $j=1,...r+2$ that $\|g_j(z)\|\le O(e(z))$.
\end{prop}
\textbf{Proof.}
1). Suppose condition (i) holds. We prove the following stronger claim by induction on $j$:
\begin{equation}\label{eq:induc_condition}
    h_{j,i}(z) \le O(e(z)) \text{ and } \nabla^k h_{j,i}(z) \le O(1), \text{ for } 1\le k\le r+2-j \ .
\end{equation}
Notice that $h_{1,i}(z)=f_i(z)$ for $i=0,...,r$ and $h_{1,i}(z)=0$ for $i\ge r+1$, thus \eqref{eq:induc_condition} holds for $j=1$. Now suppose \eqref{eq:induc_condition} holds for $j=q$. It follows from \eqref{eq:recursive_compute_h} that $$\|h_{q+1,i}(z)\|\le \sum_{l=0}^i \|\nabla h_{q,l}(z)\| \|h_{1, i-l}(z)\|\le O(e(z))\ .$$ Furthermore, for $1\le k\le r+2-q$, it follows from \eqref{eq:recursive_compute_h} and product rule of  derivative that $\nabla^k h_{q+1,i}(z)$ is a finite sum of
product of at most $(k+1)$-th order derivative of $h_{q,l}(z)$ and at most $k$-th order derivative of $h_{1, i-l}(z)$ for $l=0,\ldots, i$, all of which are $O(1)$ by induction, thus $\|\nabla^k h_{q+1,i}(z)\|\le O(1)$. This proves \eqref{eq:induc_condition} holds for $q+1$, thereby \eqref{eq:induc_condition} holds for any $q>0$ by induction. Furthermore, it follows directly from \eqref{eq:recursive} that $\|g_j(z)\|\le O(e(z))$.

2). Suppose condition (ii) holds. We show the following claims hold for any $i+j=1,...r+1$ by induction on $i+j$:
\begin{equation}\label{eq:induc_condition2}
    \|h_{j,i}(z)\| \le O(e(z)) \text{ and } \|\nabla^k h_{j,i}(z)\| \le O(1), \text{ for } 1\le k\le 2r+3-j-i \ ,
\end{equation}
then condition (i) holds by noticing $f_i(z)=h_{1,i}(z)$.
Recall that $h_{j,i}$ is recursively defined by \eqref{eq:recursive_compute_h}\eqref{eq:recursive} as shown in Figure \ref{fig:logic_flow}. For $i+j=1$, we have $h_{1,0}(z)=g_1(z)$ thus \eqref{eq:induc_condition2} holds. Now suppose \eqref{eq:induc_condition2} holds for $i+j\le q$, and we will show \eqref{eq:induc_condition2} holds for $i+j=q+1$. First, it follows from the same argument as in 1). that $h_{j,q+1-j}(z)$ for $j\ge 2$ satisfies \eqref{eq:induc_condition2} by utilizing the recursive rule \eqref{eq:recursive_compute_h}.
Now we consider the case when $j=1$. For $q\le r$, it follows from \eqref{eq:recursive} that
$$
h_{1,q}(z)= g_{q+1}(z)-\sum_{l=2}^{q+1}  \frac{1}{l!} h_{l,q}(z) \ .
$$
By utilizing the condition of $g_{q+1}(z)$ and the fact that $h_{l,q}(z)$ satisfies \eqref{eq:induc_condition2}, it holds that $h_{1,q}(z)$ satisfies \eqref{eq:induc_condition2} for $q\le r+1$. This shows condition (i) holds, thereby finishes the proof by utilizing 1). \qed

\textbf{Proof of Theorem \ref{thm:gradient-based-ODE}.} Denote $g^{(r+1)}(z):=\sum_{j=0}^{r+1} g_j(z) s^j$ as the $(r+1)$-th order Taylor series of $g(z,s)$ as in \eqref{eq:expan_g}. Then it follows from \eqref{eq:gradient-based} and $z,g(z,t)\in S$ that $$e(g(z,t))-e(z)\le \gamma \|g(z,t)-z\|\le a \gamma te(z)  \ ,$$ thus
\begin{equation}\label{eq:Fz+}
    e(g(z,t))\le (1+a \gamma t)e(z) \ .
\end{equation}
 Moreover, it follows from Proposition \ref{prop:bound_on_h} that there exists constant $c_1$ such that $\|g_j(z')\|\le c_1 e(z')$, and $\|h_{j,i}(z')\|\le c_2 e(z')$ for any $i,j$ and $z'\in \{\tz|e(\tz)\le(1+a\gamma s)e(z)\}$.


It follows from Taylor expansion of $g(z,s)$ with integral reminder that
\begin{align}\label{eq:gap1}
\begin{split}
    \|z^+ - g^{(r+1)}(z)\| = & \left\|\int_{0}^s \left.\tfrac{\partial^{r+2}}{\partial s^{r+2}} g(z, s)\right|_{s=t} \frac{t^{r+1}}{(t+1)!} dt\right\| \\
    \le & \int_{0}^s \left\|\left.\tfrac{\partial^{r+2}}{\partial s^{r+2}} g(z, s)\right|_{s=t}\right\| \frac{t^{r+1}}{(t+1)!} dt \\
    = & (r+2) \int_{0}^s \left\| g_{r+2}(g(z,t))\right\| {t^{r+1}} dt \\
    \le & c_1(r+2) \int_{0}^s e(g(z,t)) {t^{r+1}} dt \\
    \le & c_1(r+2)   e(z) \pran{\int_{0}^s{t^{r+1}} dt + a \gamma \int_{0}^s{t^{r+2}} dt }\\
    \le &c_1   e(z) \pran{s^{r+2} + a \gamma s^{r+3} } \ ,\\
\end{split}
\end{align}
where the second equality is from the definition of $g_{r+2}$, the second inequality utilizes Proposition \ref{prop:bound_on_h}, the third inequality utilizes \eqref{eq:Fz+}.

On the other hand, it follows from Taylor expansion of $Z(s)$ with integral reminder that
\begin{align}\label{eq:gap2}
\begin{split}
    &\|Z(s) - g^{(r+1)}(z)\| \\ 
    = & \norm{\sum_{j=0}^{r+1} \frac{1}{j!} \frac{d^j}{dt^j} Z(0) s^j + \int_{0}^s \frac{d^{r+2}}{dt^{r+2}} Z(t)\frac{t^{r+1}}{(r+1)!} dt - g^{(r+1)}(z)} \\
    = & \norm{\sum_{j=0}^{r+1}  \frac{1}{j!} s^j \sum_{i=0}^{rj} h_{j,i} (Z(0)) s^i + \int_{0}^s \sum_{i=0}^{r(r+2)} h_{r+2, i}(Z(t)) \frac{t^{r+1}}{(r+1)!} dt - g^{(r+1)}(z)} \\
    = & \norm{\sum_{k=0}^{r^2+2r+1} \pran{\sum_{j=0}^{\min\{k, r+1\}} \frac{1}{j!} h_{j,k-j}(z)} s^k + \int_{0}^s \sum_{i=0}^{r(r+2)} h_{r+2, i}(Z(t)) \frac{t^{r+1}}{(r+1)!} dt - g^{(r+1)}(z)} \\
    = & \norm{\sum_{k=r+2}^{r^2+2r+1} \pran{\sum_{j=0}^{\min\{k, r+1\}} \frac{1}{j!} h_{j,k-j}(z)} s^k + \int_{0}^s \sum_{i=0}^{r(r+2)} h_{r+2, i}(Z(t)) \frac{t^{r+1}}{(r+1)!} dt } \\
    \le & \sum_{k=r+2}^{r^2+2r+1} \pran{\sum_{j=0}^{r+1} \frac{1}{j!} } c_2 s^{k-r-2}  s^{r+2} e(z) + ((r+2)r +1) c_2 \int_{0}^s  \sum_{i=0}^{r(r+2)} e(Z(t)) \frac{t^{r+1}}{(r+1)!} dt  \\    
    \le & \pran{\sum_{k=0}^{r^2+r-1}  s^{k-r-2} } e c_2  s^{r+2} e(z) + ((r+2)r +1)c_2 e(z) \int_{0}^s    \frac{t^{r+1}}{(r+1)!} dt  \\    
    \le & \frac{1}{1-s} e c_2  s^{r+2} e(z) + \frac{1}{r!}   c_2 s^{r+2}e(z) \ ,   \\   
\end{split}
\end{align}
where the second equality comes from \eqref{eq:djdt}, the fourth equality utilizes the construction of $O(s^r)$ resolution ODE \eqref{eq:kth_tempo}, 
the first inequality is from Proposition \ref{prop:bound_on_h}, and the second inequity utilizes $\sum_{j=0}^k \frac{1}{j!}\le e$ and $e(Z(t))\le e(z)$ due to the
$O(s^r)$-linear-convergence condition \eqref{eq:linear-decay}.

We finish the proof by combining \eqref{eq:gap1} and \eqref{eq:gap2}. \qed




\begin{rem}\label{rem:step-size}
As shown in the proof of Theorem \ref{thm:gradient-based-ODE}, the value $c$ is upper-bounded by a polynomial of $s$ and $\|F(z)\|$, or in other words, $\|Z(s)-z^+\|$ is upper-bounded by a polynomial of $s$ and $\|F(z)\|$ with the leading term being $s^{r+2}\|F(z)\|$. Since we focus on the case when $\|F(z)\|$ is upper-bounded (by $\|F(z^0)\|$) and $s$ is small enough, the coefficient of $s^{r+2}\|F(z)\|$ in the upper-bound of $c$ dominants the condition \eqref{eq:step-size-condition}. Following a more careful calculation in Theorem \ref{thm:gradient-based-ODE}, we can obtain that the coefficient of the leading term $s^{r+2}\|F(z)\|$ in the polynomial is $O(\gamma^2)$ for both EGM and PPM. Therefore, Theorem \ref{thm:gradient-based-ODE-bound} guarantees the linear convergence rate of EGM and PPM when $\rho(s)\ge O(s^2 \gamma^3)$.
\end{rem}

Finally, we comment that the machinery stated in this section can be applied to many other algorithms for minimax problems, including but not limited to, AGDA, PDHG~\cite{chambolle2011first} and ADMM~\cite{douglas1956numerical,eckstein1992douglas}.

\section{Linear Convergence of PPM and EGM from a Discrete-Time Perspective}\label{sec:convergence}
In Corollary \ref{cor:ODE-to-DTA} (ii), we show that PPM and EGM converge linearly to a stationary solution under the $O(s)$-linear-convergence condition \eqref{eq:O(s)-strong} from a continuous-time perspective. A natural question is whether we can obtain such results within the discrete-time space. In this section, we show that a slightly modified version of the $O(s)$-linear-convergence condition can guarantee the linear convergence of PPM and EGM. The proofs completely stay in discrete-time space, although it is inspired by the convergence of their $O(s)$-resolution ODE~\eqref{eq:LypOs}. Moreover, such analysis may result in larger step-size (i.e., $s\le \frac{1}{\gamma}$ compared to $\rho(s)\ge O(s^2 \gamma^3)$ stated in Remark \ref{rem:step-size}) and does not require the high-order continuity conditions as stated in Corollary \ref{cor:ODE-to-DTA} (ii). In contrast to Corollary \ref{cor:ODE-to-DTA}, this analysis only works for convex-concave minimax problems. Similar analysis has the potential to apply to other algorithms.

\subsection{Main Results}

First, we define a variant of $O(s)$-linear-convergence condition \eqref{eq:O(s)-strong}:

\begin{mydef}
Define $\FF=\left\{F(z_1)+F(z_2)|z_1,z_2\in\RRmn\right\}$. We say the minimax function $L(x,y)$ 
satisfies the strong $O(s)$-linear-convergence condition of PPM and EGM if there exists $\rho(s)>0$ such that it holds for any $c\in \FF$ and $z=(x,y)\in \RRmn$ that
\begin{equation}\label{eq:O(s)-strong-PPM}
    c^T \left[\begin{matrix} A(z)- \frac{s}{2} A(z)^2+ \frac{s}{2} B(z)B(z)^T & 0 \\ 0 & C(z)- \frac{s}{2} C(z)^2+ \frac{s}{2} B(z)^T B(z) \end{matrix}\right] c \ge \frac{1}{2} \rho(s) \|c\|^2 \ .
\end{equation}
\end{mydef}

Compared to \eqref{eq:OsconditionPPM}, \eqref{eq:O(s)-strong-PPM} is a slightly stronger condition in the sense that $c$ in \eqref{eq:O(s)-strong-PPM} is chosen from a larger space $\FF$ compared with that in \eqref{eq:OsconditionPPM}. 

Theorem \ref{thm:ppm} presents the linear convergence rate of PPM \eqref{eq:ppm} when the function $L(x,y)$ satisfies the strong $O(s)$-linear-convergence condition \eqref{eq:O(s)-strong-PPM}.
\begin{thm}\label{thm:ppm}
\textbf{(Fast convergence of PPM)} Consider PPM with iterate update \eqref{eq:ppm} and step-size $s\le \frac{1}{3\gamma}$. Suppose $L(x,y)$ is convex-concave and it satisfies the strong $O(s)$-linear-convergence condition \eqref{eq:O(s)-strong-PPM}, then it holds for all iteration $k\ge 0$ that
\begin{equation*}
        \normt{\Fzk} \le \pran{\frac{1-\frac{s\rho(s)}{2}}{1+\frac{s\rho(s)}{4}}}^k \normt{F(z_0)} \ .
\end{equation*}
\end{thm}

\begin{rem}
Theorem \ref{thm:ppm} shows that PPM with step-size $s\le \frac{1}{3\gamma}$ finds a solution $z$ such that $\|F(z)\|^2\le \varepsilon$ within $O(\frac{1}{s\rho(s)}\log(\frac{1}{\varepsilon}))$ iterations. 
\end{rem}

Now we turn to EGM. Our first result is Theorem \ref{thm:slow-mp}, which shows that when the step-size is small enough such that $\rho(s)\ge  8 s^2\gamma^3$, EGM has linear convergence.
The linear convergence rate is slower than that of PPM (Theorem \ref{thm:ppm}) due to the required  smaller step-size to satisfy  $\rho(s)\ge  8 s^2\gamma^3$. Secondly, in the case when the $L(x,y)$ is a convex-concave quadratic function, Theorem \ref{thm:fast-mp} shows that EGM can take a larger step-size, which recovers the same order of linear convergence rate of PPM in Theorem \ref{thm:ppm}. We further compare the slow rate and fast rate in Remark \ref{rem:sl_vs_fa}.
\begin{thm}\label{thm:slow-mp}
\textbf{(Slow convergence of EGM)} Consider the EGM with iterate update \eqref{eq:mp} and step-size $s$. Suppose $L(x,y)$ is convex-concave and it satisfies the $O(s)$-linear-convergence condition \eqref{eq:O(s)-strong-PPM}, and suppose the step-size $s$ satisfies $s\le\frac{1}{2\gamma}$ and $\rho(s)\ge  8 s^2\gamma^3$, then it holds for all iteration $k\ge 0$ that
\begin{equation*}
        \normt{\Fzk}\le \pran{\frac{1-\frac{s\rho(s)}{5}}{1+\frac{s\rho(s)}{5}}}^k \normt{F(z_0)} \ .
\end{equation*}
\end{thm}

\begin{thm}\label{thm:fast-mp}
\textbf{(Fast convergence of EGM for quadratic function)} Consider the EGM with iterate update \eqref{eq:mp} and step-size $s$. Suppose $L(x,y)$ is a quadratic function
\begin{equation}\label{eq:quadratic}
    L(x,y)=\frac{1}{2}x^T A x +x^T B y -\frac{1}{2} y^T C y + d^T x + e^T y\ ,
\end{equation}
where matrix $A$ and $C$ are positive semi-definite matrices. Suppose $L(x,y)$ satisfies the $O(s)$-linear-convergence condition \eqref{eq:O(s)-strong}, and suppose the step-size $s$ satisfies $s\le\frac{1}{8\gamma}$, then it holds for all iteration $k\ge 0$ that
\begin{equation*}
        \normt{\Fzk}\le \pran{\frac{1-\frac{s\rho(s)}{5}}{1+\frac{s\rho(s)}{5}}}^k \normt{F(z_0)} \ .
\end{equation*}
\end{thm}

\begin{rem}\label{rem:sl_vs_fa} 
Here we compare the slow rate (Theorem \ref{thm:slow-mp}) and fast rate (Theorem \ref{thm:fast-mp}) of EGM. Recall that Theorem \ref{thm:fast-mp} (fast rate) requires $s\le \frac{1}{8\gamma}$, while Theorem \ref{thm:slow-mp} and Remark \ref{rem:step-size} (slow rate) requires
\begin{equation}\label{eq:weak_condition}
    \rho(s)\ge 8s^2\gamma^3 \ .
\end{equation}
Let us consider the two standard scenarios discussed in the introduction section. When $L(x,y)$ is $\mu$-strongly convex-strongly concave, $\rho(s)\ge \mu$, condition \eqref{eq:weak_condition} requires that $s\sim O\pran{\sqrt{\frac{\mu}{\gamma^3}}}$, thus to find a solution $z$ such that $\|F(z)\|^2\le \varepsilon$, Theorem \ref{thm:slow-mp} suggests EGM needs $O\pran{\pran{\frac{\gamma}{\mu}}^{3/2}\log\pran{\frac{1}{\varepsilon}}}$ iterations. In contrast, Theorem \ref{thm:fast-mp} suggests EGM needs $\pran{\frac{\gamma}{\mu}}\log\pran{\frac{1}{\varepsilon}}$ iterations. When $L(x,y)=y^T B x$, $\rho(s)=\lambda_{\min}^+ (BB^T)$, condition \eqref{eq:weak_condition} requires that $s\sim O\pran{\frac{\lambda_{\min}^+ (BB^T)}{\gamma^3}}$, thus to find a solution $z$ such that $\|F(z)\|^2\le \varepsilon$, Theorem \ref{thm:slow-mp} suggests EGM needs $O\pran{\pran{\frac{\gamma^2}{\lambda_{\min}^+ (BB^T)}}^{3}\log\pran{\frac{1}{\varepsilon}}}$ iterations. In contrast, Theorem \ref{thm:fast-mp} suggests EGM needs $O\pran{\pran{\frac{\gamma^2}{\lambda_{\min}^+ (BB^T)}}\log\pran{\frac{1}{\varepsilon}}}$ iterations. \red{Finally, we comment that the different step-size requirement for the quadratic and the general objective is also observed for GDA to solve one-side strongly convex minimax problems~\cite{du2017stochastic,du2019linear}.}
\end{rem}

\red{
\subsection{Proof Scratch of Theoerm \ref{thm:ppm}-\ref{thm:fast-mp}}\label{sec:scratch-proof}
Here we provide a proof scratch of the linear convergence of PPM and EGM (Theoerm \ref{thm:ppm}-\ref{thm:fast-mp}). The proofs for these three theorems have very similar structures and they are all inspired by the energy decay of their $O(s)$-resolution ODE \eqref{eq:LypOs}. 

We consider the discrete-time counterpart of the energy function~\eqref{eq:continuous-energy} and studies its decay in discrete-time under their $O(s)$-linear-convergence-conditions. Notice that
\begin{align*}
    \frac{1}{2}\|F(z_{k+1})\|^2 - \frac{1}{2}\|F(z_{k})\|^2 = \frac{1}{2}\pran{\sumF}^T \pran{\minusF} \ .
\end{align*}

The first step in the proof is to show that there exists $R(z_k,s)\in \RR^{(n+m)\times(n+m)}$ such that
\begin{align*}
    \minusF =  s R(z_k,s) \pran{\sumF} \ .
\end{align*}
Now suppose $R(z_k,s)$ has Taylor expansion of $s$: $R(z_k,s)=\sum_{j=0}^\infty R_j(z_k) s^j$. 
Indeed, it turns out the first two terms in the Taylor expansion of $R(z_k,s)$ for PPM and EGM after canceling out the skew-symmetric interaction terms is exactly the term in $O(s)$-linear-convergence condition \eqref{eq:O(s)-strong-PPM}:
\begin{equation*}
    -\left[\begin{matrix} A(z_k)- \frac{s}{2} A(z_k)^2+ \frac{s}{2} B(z_k)B(z_k)^T & 0 \\ 0 & C(z_k)- \frac{s}{2} C(z_k)^2+ \frac{s}{2} B(z_k)^T B(z_k) \end{matrix}\right] \ .
\end{equation*}
This is not surprising due to the construction of the $O(s^r)$-resolution ODE. Thus it follows \eqref{eq:O(s)-strong-PPM} that 
\begin{align}\label{eq:discrete-energy-decay}
\begin{split}
        &\frac{1}{2}\|F(z_{k+1})\|^2 - \frac{1}{2}\|F(z_{k})\|^2\\
    = &  s\pran{\sumF}^T \sum_{j=0}^\infty R_j(z_k) s^j \pran{\sumF} \\
    = & s\pran{\sumF}^T \left[\begin{matrix} A- \frac{s}{2} A^2+ \frac{s}{2} BB^T & 0 \\ 0 & C- \frac{s}{2} C^2+ \frac{s}{2} B^T B \end{matrix}\right] \pran{\sumF} \\
    & + s\pran{\sumF}^T \sum_{j=2}^\infty R_j(z_k) s^j \pran{\sumF} \\
    \le & -s\pran{1-\frac{1}{2}\rho(s)} \|\sumF\|^2 + s\pran{\sumF}^T \sum_{j=2}^\infty R_j(z_k) s^j \pran{\sumF} \ ,
\end{split}
\end{align}
where we omit $z_k$ as arguments in $A, B, C$ in the third equality for notational convenience. The first term in the right-hand side of \eqref{eq:discrete-energy-decay} provides a sufficient decay of the energy function, which results in the linear convergence of PPM/EGM. The rest of the proof is to show that the last sum term in the right-hand side of \eqref{eq:discrete-energy-decay} (i.e., the $o(s^2)$ terms in the Taylor expansion) does not affect this linear rate much (upto a constant) when the step-size is small enough.

For the slow rate (such as Theorem \ref{thm:slow-mp}), we show that $\|\sum_{j=2}^\infty R_j(z_k) s^j\|\le c_3 s^2$ for a constant $c_3$, thereby the linear rate holds as long as $\rho(s)\ge 2 c_3 s^2$. Such an argument is very general and can be applied to analyze other algorithms. This is consistent with and provides a different perspective of the linear rate stated in Corollary \ref{cor:ODE-to-DTA}.

The fast rate (such as Theorem \ref{thm:ppm} and Theorem \ref{thm:fast-mp}) allows for larger step-size ($s\le O(1/\gamma)$), but requires more subtle calculations, which may not hold in general. In this argument, we show that as long as $s\le O(1/\gamma)$, it holds that
\begin{align*}
    &\left|\pran{\sumF}^T \sum_{j=2}^\infty R_j(z_k) s^j \pran{\sumF}\right| \\
    \le & \frac{1}{2} \pran{\sumF}^T \left[\begin{matrix} A- \frac{s}{2} A^2+ \frac{s}{2} BB^T & 0 \\ 0 & C- \frac{s}{2} C^2+ \frac{s}{2} B^T B \end{matrix}\right] \pran{\sumF} \ ,
\end{align*}
thus the last term in RHS of \eqref{eq:discrete-energy-decay} results in at most a factor of $2$ in the analysis. We also want to mention that, as shown in the proof later, the analysis of the above inequality can be highly nontrivial and heavily depends on the properties of generalized block skew-symmetric matrices, which we define and explain in Appendix \ref{sec:GBSS}.

Notice that EGM has a fast rate for quadratic problems while it has a slow rate for general problems. From the proof perspective, this is because $R(z_k, s)$ has a complicated expression for general problems, which can be greatly simplified for quadratic problems. 

The formal proofs of Theorem \ref{thm:ppm}-\ref{thm:fast-mp} are left in Appendix \ref{sec:proofs}.
}




\section{Conclusion and Future Directions}

In this paper, we present a new machinery -- an $O(s^r)$-resolution ODE framework -- for analyzing the behavior of a generic DTA, and apply it to unconstrained minimax problems. We propose the $r$-th degree ODE expansion of a DTA to construct the unique $O(s^r)$-resolution ODE. From the $O(s^r)$-resolution ODE, we present how to obtain an $O(s^r)$-linear-convergence condition with respect to an energy function, which not only guarantees the linear convergence of the $O(s^r)$-resolution ODE, but also guarantees the linear convergence of the original DTA if the energy function is chosen properly. We utilize this machinery to study GDA, PPM and EGM for solving minimax problems, which provides intuitive explanations of their different behaviors and also results in tighter conditions under which these methods have linear convergence. This machinery can also help design new algorithms. 

Future directions of this line of research include (i) using this machinery to study other algorithms, for example, PDHG, ADMM, etc; (ii) extending this machinery to other settings, for example, constrained optimization and stochastic algorithms; (iii) extending this machinery to bridge the sublinear convergence of a DTA and its corresponding ODEs. Furthermore, we present Conjecture \ref{con:inf}. Suppose it is true, then we can utilize an ODE to fully represent a DTA.


\section*{Acknowledgement}
The author would like to express his gratitude to Robert M. Freund for reading an early version of the paper and for thoughtful discussions that helped to position the paper. The author also wishes to thank Renbo Zhao, Ben Grimmer, Miles Lubin, Oliver Hinder and David Applegate for helpful discussions. The author would like to thank the anonymous referees and the associate editor for the constructive feedback, which results in a significantly improved version of the manuscript.

\bibliographystyle{amsplain}
\bibliography{LF-papers}

\appendix
\section{Appendix}




\subsection{$O(s)$-Linear-Convergence Condition of $L(x,y)= f(C_1 x) + x^T B y - g(C_2 y)$}\label{sec:app-example}
\begin{prop}

Consider $L(x,y)= f(C_1 x) + x^T B y - g(C_2 y)$. Define
$$
a_1=\left\{\begin{array}{cl}
     \min\pran{\mu \lambda_{\min}^+ (C_1^T C_1), s\lambda_{\min}^+ (B B^T)} & \ \emph{if} \ \sin\pran{\eRange{B}, \eRange{C_1^T}}= 0  \\
     \min\pran{\mu \lambda_{\min}^+ (C_1^T C_1)\sin^2\pran{\eRange{B}, \eRange{C_1^T}}, s\lambda_{\min}^+ (B B^T)}& \ \emph{otherwise} \ ,
\end{array} \right.
$$
and
$$
a_2=\left\{\begin{array}{cl}
     \min\pran{\mu \lambda_{\min}^+ (C_2^T C_2), s\lambda_{\min}^+ (B^T B)} & \ \emph{if} \ \sin\pran{\eRange{B^T}, \eRange{C_2^T}}= 0  \\
     \min\pran{\mu \lambda_{\min}^+ (C_2^T C_2)\sin^2\pran{\eRange{B^T}, \eRange{C_2^T}}, s\lambda_{\min}^+ (B^T B)}& \ \emph{otherwise}\  \ ,
\end{array} \right.
$$
where $\sin(\cdot, \cdot)$ is the cosine angle between two linear spaces\footnote{Suppose $\mathcal{A},\mathcal{B}$ are two linear subspaces in $\RR^m$, then $\cos(\mathcal{A},\mathcal{B}):=\min_{a\in\mathcal{A}, b\in\mathcal{B}} \cos(a,b)$, and $\sin(\mathcal{A},\mathcal{B})=\sqrt{1-\cos^2(\mathcal{A},\mathcal{B})}$.}. Then $L(x,y)$ satisfies the $O(s)$-linear-convergence condition with $\rho(s) \ge \min\{a_1,a_2\}>0$.
\end{prop}

{\bf Proof.} 
Suppose it holds for any $x\in \Range{C_1^T} + \Range{B}$ that 
\begin{equation}\label{eq:suppose_hold}
    x^T \pran{\nabla_{xx}L(x,y)+ s \nabla_{xy}L(x,y)\nabla_{xy}L(x,y)^T} x \ge a_1 \|x\|^2 \ ,
\end{equation}
then symmetrically for any $y\in \Range{C_2^T} + \Range{B^T}$ it holds that
\begin{equation*}
    y^T \pran{\nabla_{yy}L(x,y)+ s \nabla_{xy}L(x,y)^T\nabla_{xy}L(x,y)} y \ge a_2\|y\|^2 \ ,
\end{equation*}
which proves \eqref{eq:O(s)-strong} with $\rho(s)=\min\{a_1,a_2\}>0$ by noticing $\FF\subseteq \pran{\Range{C_1^T} + \Range{B}}\times \pran{\Range{C_2^T} + \Range{B^T}}$. 
Now let us prove \eqref{eq:suppose_hold}. First, notice that $\nabla_{xx}L(x,y)\succeq \mu C_1^T C_1$ and $\nabla_{xy}L(x,y)=B$, thus we just need to show
\begin{equation}\label{eq:suppose_hold2}
    x^T \pran{\mu C_1^T C_1+ s BB^T} x \ge a_1 \|x\|^2 \ .
\end{equation}
If $\sin\left(B B^T, \Cop\right)=0$, then either $x\in \Range{C_1^T}$ thus $x^T \pran{\mu C_1^T C_1+ s BB^T} x \ge \mu \lambda_{\min}^+ (C_1^T C_1)\|x\|^2$, or $x\in \Range{B}$ thus $x^T \pran{\mu C_1^T C_1+ s BB^T} x \ge s \lambda_{\min}^+ (BB^T)\|x\|^2$. In either case \eqref{eq:suppose_hold2} holds.

If $\sin\left(B B^T, \Cop\right)\not=0$, suppose $x=x_1+x_2$ where $x_1\in \Range{B^T}$ and $x_2\in \Range{C_1^T}$. It is obvious that \eqref{eq:suppose_hold2} holds if $x_2 = 0$. Now define $\PB(x)=B(BB^T)^{+}B^Tx$ as the projection operator onto $\Range{B}$, and $\PB^T(x)=x-\PB(x)$ be the projection operator onto the perpendicular space of $\Range{B}$, then it holds that
\begin{align*}
    \begin{split}
        &x^T \pran{\mu C_1^T C_1+ s BB^T} x \\
        = & (x_1+\PB(x_2)+\PB^T(x_2))^T \pran{\mu C_1^T C_1+ s BB^T} (x_1+\PB(x_2)+\PB^T(x_2)) \\
        = & (x_1+\PB(x_2))^T \pran{\mu C_1^T C_1+ s BB^T} (x_1+\PB(x_2)) + \mu (\PB^T(x_2))^T C_1^T C_1 \PB^T(x_2) \\
        \ge & (x_1+\PB(x_2))^T \pran{s BB^T} (x_1+\PB(x_2)) + \mu (\PC(\PB^T(x_2)))^T C_1^T C_1 \PC(\PB^T(x_2)) \\
        \ge & s\lambda_{\min}^+(BB^T)\|x_1+\PB(x_2)\|^2 + \mu \lambda_{\min}^+(C_1^T C_1) \|(\PC(\PB^T(x_2)))\|^2 \\
        \ge & a_1\|x_1+\PB(x_2)\|^2 + \mu \lambda_{\min}^+(C_1^T C_1) \sin^2\pran{\Range{B}, \Range{C_1^T}} \|\PB^T(x_2)\|^2 \\
        \ge & a_1\|x_1+\PB(x_2)\|^2 + a_1 \|\PB^T(x_2)\|^2 \\
        = & a_1 \|x\|^2\ ,
    \end{split}
\end{align*}
where the second equality uses $B^T \PB^T(x_2)=0$, the first inequality is from $(x_1+\PB(x_2))^T \pran{\mu C_1^T C_1} (x_1+\PB(x_2))\ge 0$ and $C_1 \PC^T(\PB^T(x_2))=0$, the second inequality is because $x_1+\PB(x_2)\in \Range{B^T}$  and $\PC(\PB^T(x_2))\in\Range{C_1^T}$, the third inequality uses the definition of $a_1$ and the definition of $\cos$ between two space, the fourth inequality is due to the definition of $a_1$, and the last equality is from $x_1+\PB(x_2)\in \Range{B^T}$ and  $\PB^T(x_2)\perp \Range{B^T}$. This finishes the proof. \qed

\section{Taylor Expansion of Operator $(I+sF)^{-1}$}\label{sec:expansion_operator}
\red{
Here we derive the third order Taylor expansion of operator $(I+sF)^{-1}$ as stated in \eqref{eq:second-order-Taylor-PPM}.
Suppose $(I+sF)^{-1}(z)=g_0(z)+g_1(z)s+g_2(z)s^2+g_3(z) s^3 + o(s^3)$, then it holds that
\begin{align}\label{eq:inverse_expansion}
    \begin{split}
        z&=(I+sF)(g_0(z)+g_1(z)s+g_2(z)s^2+g_3(z) s^3) + o(s^3)\\
        &= g_0(z)+g_1(z)s+g_2(z)s^2+g_3(z) s^3+ sF(g_0(z)+g_1(z)s+g_2(z)s^2)+ o(s^3)\ .
    \end{split}
\end{align}
By comparing the $O(1)$ term in both sides of \eqref{eq:inverse_expansion}, we have
$g_0(z)=z$. By comparing the $O(s)$ term in both sides of \eqref{eq:inverse_expansion}, we have
$$0=g_1(z)+F(g_0(z))=g_1(z)+F(z)\ ,$$thus $g_1(z)=-F(z)$. Notice $F(g_0(z)+sg_1(z))=F(z-sF(z))=F(z)-s\nabla F(z) F(z)+o(s)$. By comparing the $O(s^2)$ term in both side of \eqref{eq:inverse_expansion}, we have
$$0=g_2(z)-\nabla F(z) F(z)\ , $$ thus $g_2(z)=\nabla F(z) F(z)$. Notice
\begin{align*}
    \begin{split}
       &F(g_0(z)+g_1(z)s+ g_2(z)s^2)\\
       =&F(z-sF(z)+s^2 \nabla F(z) F(z))\\
       =&F(z) +\nabla F(z)( -sF(z)+s^2\nabla F(z) F(z)) + \frac{1}{2}\nabla^2 F(z) (sF(z), sF(z)) + o(s^2) \\
       =&F(z)-s\nabla F(z) F(z) + s^2 \pran{(\nabla F(z))^2 F(z) + \frac{1}{2}\nabla^2 F(z) (F(z), F(z))} + o(s^2) \ .
    \end{split}
\end{align*}
By comparing the $O(s^3)$ term in both sides of \eqref{eq:inverse_expansion}, we have $$0=g_3(z)+(\nabla F(z))^2 F(z) +\frac{1}{2} \nabla^2 F(z) (F(z), F(z)) \ ,$$ thus $g_3(z)=-(\nabla F(z))^2 F(z) -\frac{1}{2} \nabla^2 F(z) (F(z), F(z))$, which yield \eqref{eq:second-order-Taylor-PPM}.
}

\section{Generalized Block Skew-Symmetric Matrix and Its Basic Properties}\label{sec:GBSS}

Here is the definition of a generalized block skew-symmetric matrix:
\begin{mydef}
We say a matrix $M\in\RR^{(n+m)\times (n+m)}$ is generalized block skew-symmetric if $M$ has the structure: $M=\left[\begin{matrix}  A & B\\-B^T  & C \end{matrix}\right]$ where $A\in \RR^{n\times n}, C\in \RR^{m\times m}$ are symmetric matrices and $B\in \RR^{n\times m}$ is an arbitrary matrix.
\end{mydef}

\begin{rem}
Going back to the minimax problem, $\nabla F(z)=\left[\begin{matrix}  \nabla_{xx}L(x,y) & \nabla_{xy}L(x,y)\\ -\nabla_{xy}L(x,y)^T & \nabla_{yy}L(x,y) \end{matrix}\right]$ is a generalized block skew-symmetric matrix for any $z$. 
\end{rem}

\vspace{0.2cm}

Let $M=\ABC$ be a generalized symmetric matrix. Denote $M^i=\left[\begin{matrix}\Mp{i}{11} & \Mp{i}{12}\\  \Mp{i}{21}  & \Mp{i}{22} \end{matrix}\right]$ as the $i$th power of matrix $M$, where $\Mp{i}{jl}$ for $j,l\in \{1,2\}$ is the corresponding block of $M^i$. In particular, we define $M^0$ to be the identity matrix. The next proposition shows that $M^i$ keeps the generalized block skew-symmetricity.

\begin{prop}\label{prop:power-symmetry}
Suppose $M$ is a generalized block skew-symmetric matrix, then for any positive integer $i$, $M^i$ is a generalized block skew-symmetric matrix.
\end{prop}
{\bf Proof.} We'll prove the Proposition \ref{prop:power-symmetry} by induction. First notice that Proposition \ref{prop:power-symmetry} is satisfied with $i=1$. Now suppose Proposition \ref{prop:power-symmetry} is satisfied with $i$. Notice that
\begin{equation}\label{eq:silly}
M^{i+1}=M M^{i} =M^{i} M \ ,    
\end{equation}
which yield the following update by matrix multiplication rules:
\begin{align}\label{eq:recursive_update}
\begin{split}
     \Mp{i+1}{11}&=A\Mp{i}{11}+B\Mp{i}{21}=\Mp{i}{11}A-\Mp{i}{12}B^T, \\
 \Mp{i+1}{12}&=A\Mp{i}{12}+B\Mp{i}{22}=\Mp{i}{11}B+\Mp{i}{12}C,\\
 \Mp{i+1}{21}&=-B^T\Mp{i}{11}+C\Mp{i}{21}=\Mp{i}{21}A-\Mp{i}{22}B^T ,\\ 
 \Mp{i+1}{22}&=-B^T\Mp{i}{12}+C\Mp{i}{22}=\Mp{i}{21}B+\Mp{i}{22}C.
\end{split}
\end{align}
Therefore,
$$\Mp{i+1}{11}=\frac{1}{2}\left(A\Mp{i}{11}+B\Mp{i}{21}+\Mp{i}{11}A-\Mp{i}{12}B^T\right)= \frac{1}{2}\left(\pran{A\Mp{i}{11}+B\Mp{i}{21}}+ \pran{A\Mp{i}{11}+B\Mp{i}{21}}^T\right)$$
is symmetric. Similarly, we have $\Mp{i+1}{22}$ is symmetric. Meanwhile, it holds that
$$
\Mp{i+1}{12}=A\Mp{i}{12}+B\Mp{i}{22}=-\left(\Mp{i}{21}A-\Mp{i}{22}B^T\right)^T = -\pran{\Mp{i+1}{21}}^T \ ,
$$
which finishes the proof for (i) by induction. \qed

\vspace{0.2cm}
















The next proposition provides upper and lower bounds on $\Mp{i}{11}$ and $\Mp{i}{22}$:
\begin{prop}\label{prop:power_bound}
Suppose $M$ is a generalized block skew-symmetric matrix, and $\|M\|\le\gamma$, then it holds for $i\ge 3$ that
\begin{equation}\label{eq:lem_bound_1}
    -(i-1)\gamma^{i-2}(\gamma A+ BB^T) \preceq \Mp{i}{11}\preceq (i-1)\gamma^{i-2}(\gamma A+ BB^T) \ ,
\end{equation}
and
\begin{equation}\label{eq:lem_bound_2}
    -(i-1)\gamma^{i-2}(\gamma C+ B^TB) \preceq \Mp{i}{22}\preceq (i-1)\gamma^{i-2}(\gamma C+ B^TB) \ .
\end{equation}
Furthermore, it holds for any integer $i\ge 3$ and $c\in\RRmn$ that
\begin{equation*}
    \pabs{c^T M^i c} \le (i-1)\gamma^{i-2} c^T \left[\begin{matrix}  \gamma A+ BB^T & 0\\0  & \gamma C+ B^T B \end{matrix}\right]c \ .
\end{equation*}
\end{prop}

The following two facts will be needed for the proof of Proposition \ref{prop:power_bound}.

\begin{fact}\label{fact:average}
    Suppose $S_1$ and $S_2$ are symmetric matrices, then
    $$
    -(S_1^2 + S_2^2)\preceq S_1S_2+S_2S_1 \preceq S_1^2 + S_2^2 \ .
    $$
\end{fact}
{\bf Proof.}
It is easy to check that
$$
S_1^2 + S_2^2 - S_1S_2+S_2S_1 = (S_1-S_2)^T(S_1-S_2) \succeq 0 ,
$$
and
$$
S_1^2 + S_2^2 + S_1S_2+S_2S_1 = (S_1+S_2)^T(S_1+S_2) \succeq 0 ,
$$
which finishes the proof by rearranging the above two matrix inequalities. \qed

\begin{fact}\label{fact:Mioo}
Suppose $M$ is a generalized block skew-symmetric matrix, then
\begin{equation}\label{eq:fact2}
    \Mp{i}{11}=A\Mp{i-2}{11}A-B\Mp{i-2}{22} B^T  - \pran{\sum_{j=0}^{i-3} B\Mp{j}{22}B^T A^{i-2-j} + A^{i-2-j} B\Mp{j}{22}B^T}\ .
\end{equation}
\end{fact}
{\bf Proof.}
By recursively using the update rule \eqref{eq:recursive_update} and rearranging the equality, it holds that:
\begin{align*}
    \begin{split}
        \Mp{i}{11} & = A \Mp{i-1}{11} + B\Mp{i-1}{21} \\
        & = A(\Mp{i-2}{11}A-A\Mp{i-2}{12}B^T) + B(\Mp{i-2}{21}A-\Mp{i-2}{22}B^T) \\
        & = A\Mp{i-2}{11}A-B\Mp{i-2}{22} B^T + \pran{B\Mp{i-2}{21}A - A\Mp{i-2}{12}B^T} \\
        & = A\Mp{i-2}{11}A-B\Mp{i-2}{22} B^T + \pran{B\Mp{i-3}{21}A^2 - A^2\Mp{i-3}{12}B^T} - \pran{B\Mp{i-3}{22}B^T A + A B\Mp{i-3}{22}B^T} \\
        & = \cdots \\
        & = A\Mp{i-2}{11}A-B\Mp{i-2}{22} B^T + \pran{BB^T A^{i-2} + A^{i-2} BB^T} - \pran{\sum_{j=1}^{i-3} B\Mp{j}{22}B^T A^{i-2-j} + A^{i-2-j} B\Mp{j}{22}B^T} \\
        & = A\Mp{i-2}{11}A-B\Mp{i-2}{22} B^T  - \pran{\sum_{j=0}^{i-3} B\Mp{j}{22}B^T A^{i-2-j} + A^{i-2-j} B\Mp{j}{22}B^T} \ .
    \end{split}{}
\end{align*}{}\qed

Now let us go back to the proof of Proposition \ref{prop:power_bound}.

{\bf Proof of Proposition \ref{prop:power_bound}.}

Notice that $A$ is positive semi-definite and $\|M\|=\gamma$, thus $\|A\|\le \gamma$ and $\|\Mp{i-2}{11}\|\le \gamma^{i-2}$, whereby $A^{1/2}\Mp{i-2}{11}A^{1/2}\preceq \gamma^{i-1}I$. Therefore, it holds that
\begin{equation}\label{eq:temp_1}
    0\preceq \tfrac{1}{\gamma^i} A\Mp{i-2}{11}A = \tfrac{1}{\gamma^i} A^{1/2}\left(A^{1/2}\Mp{i-2}{11}A^{1/2}\right)A^{1/2} \preceq \tfrac{1}{\gamma} A\ .
\end{equation}

Notice that $\Mp{i-2}{22}\preceq \gamma^{i-2}I$, thus it holds that
\begin{equation}\label{eq:temp_2}
    0\preceq \tfrac{1}{\gamma^i} B\Mp{i-2}{22}B^T = \tfrac{1}{\gamma^i} B\Mp{i-2}{22}B^T \preceq \tfrac{1}{\gamma^{2}} BB^T\ .
\end{equation}

For any $0\le j\le i-3$, we have from Fact \ref{fact:average} by choosing $S_1=\tfrac{1}{\gamma^{2+j}} B\Mp{j}{22}B^T$ and $S_2=\tfrac{1}{\gamma^{i-j-2}} A^{i-j-2}$ that
\begin{align}\label{eq:temp_3}
    \begin{split}
        & \  \tfrac{1}{\gamma^i} B\Mp{j}{22}B^T A^{i-2-j} + \tfrac{1}{\gamma^i} A^{i-2-j} B\Mp{j}{22}B^T \\
        \preceq & \ \left(\tfrac{1}{\gamma^{2+j}} B\Mp{j}{22}B^T\right)^2 + \left( \tfrac{1}{\gamma^{i-j-2}} A^{i-j-2}\right)^2 \\
        = & \ \tfrac{1}{\gamma^{2j+4}} B\left(\Mp{j}{22}B^T B \Mp{j}{22}\right) B^T + \tfrac{1}{\gamma^{2i-2j-4}} A^{1/2} A^{2i-2j-3} A^{1/2} \\
        \preceq  & \ \tfrac{1}{\gamma^{2}} B B^T + \tfrac{1}{\gamma}  A\ ,
    \end{split}
\end{align}
where the second matrix inequality is because $B^TB\preceq\gamma^2 I$, $\Mp{j}{22}\preceq \gamma^j I$ and $A\preceq \gamma I$. Similarly, it holds that
\begin{equation}\label{eq:temp_4}
    \tfrac{1}{\gamma^i} B\Mp{j}{22}B^T A^{i-2-j} + \tfrac{1}{\gamma^i} A^{i-2-j} B\Mp{j}{22}B^T \succeq  -\tfrac{1}{\gamma^{2}} B B^T -\tfrac{1}{\gamma} A .
\end{equation}
Substituting \eqref{eq:temp_1}, \eqref{eq:temp_2}, \eqref{eq:temp_3} and \eqref{eq:temp_4} into \eqref{eq:fact2} yields
\begin{align}
    \begin{split}
        \tfrac{1}{\gamma^i} \Mp{i}{11} & =  \tfrac{1}{\gamma^i}\left(A\Mp{i-2}{11}A-B\Mp{i-2}{22} B^T  - \pran{\sum_{j=0}^{i-3} B\Mp{j}{22}B^T A^{i-2-j} + A^{i-2-j} B\Mp{j}{22}B^T}\right) \\
        &\preceq  \pran{ \tfrac{1}{\gamma} A + \tfrac{1}{\gamma^2} BB^T + (i-2)(\tfrac{1}{\gamma} A+\tfrac{1}{\gamma^{2}} BB^T)} \\
        &= (i-1)(\tfrac{1}{\gamma} A+\tfrac{1}{\gamma^{2}} BB^T) \ ,
    \end{split}
\end{align}
and
\begin{align}
    \begin{split}
        \tfrac{1}{\gamma^i} \Mp{i}{11} & =  \tfrac{1}{\gamma^i}\left(A\Mp{i-2}{11}A-B\Mp{i-2}{22} B^T  - \pran{\sum_{j=0}^{i-3} B\Mp{j}{22}B^T A^{i-2-j} + A^{i-2-j} B\Mp{j}{22}B^T}\right) \\
        &\succeq  \pran{-\tfrac{1}{\gamma} A -\tfrac{1}{\gamma} BB^T - (i-2)(\tfrac{1}{\gamma} A+\tfrac{1}{\gamma^{2}} BB^T)} \\
        &= -(i-1)(\tfrac{1}{\gamma} A+\tfrac{1}{\gamma^{2}} BB^T) \ .
    \end{split}
\end{align}
which furnishes the proof of \eqref{eq:lem_bound_1}. The proof of \eqref{eq:lem_bound_2} can be obtained symmetrically. Furthermore, it follows from Proposition \ref{prop:power-symmetry} that $M^i$ is generalized block skew-symmetric, thus
\begin{equation}
    \pabs{c^T M^i c} = \pabs{c^T \left[\begin{matrix}  \Mp{i}{11} & 0\\0  & \Mp{i}{22} \end{matrix}\right] c } \le (i-1)\gamma^{i-2} c^T \left[\begin{matrix}  \gamma A+ BB^T & 0\\0  & \gamma C+ B^T B \end{matrix}\right]c \ ,
\end{equation}
which finishes the proof of Proposition \ref{prop:power-symmetry}.

\qed

\section{Proofs in Section \ref{sec:convergence}}\label{sec:proofs}
\subsection{Proof of Theorem \ref{thm:ppm}}



The following two propositions will be needed for the proof of Theorem \ref{thm:ppm}.
\begin{prop}\label{prop:minusF}
For given $z$ and $\hz$, let $M=\int_{0}^1 \nabla F(z+t (\hz-z))dt$, then $F(\hz)-F(z)=M(\hz-z)$.
\end{prop}{}
{Proof.} Let $\phi(t)=F(z+t (\hz-z))dt$, then $\phi(0)=F(z)$, $\phi(1)=F(\hz)$ and $\phi'(t)=\nabla F(z+t(\hz-z))(\hz-z)$. From the fundamental theorem of calculus, we have
$$
F(\hz)-F(z)=\phi(1)-\phi(0)=\int_{0}^1 \phi'(t)dt = \int_{0}^1 \nabla F(z+t(\hz-z))(\hz-z)dt = M(\hz-z)\ .
$$
\qed
\begin{prop}\label{prop:sumF}
Consider PPM with iterate update \eqref{eq:ppm} and step-size $s\le \frac{1}{3\gamma}$, then for any iteration $k$, it holds that
    $$\|\sumF\|^2 \ge 2\|F(z_k)\|^2+\|F(z_{k+1})\|^2 \ .$$
\end{prop}{}
{\bf Proof.} Let $M=\int_{0}^1 \nabla F(z_{k+1}+t (z_{k+1}-z_k))dt$, then $\|M\|\le\int_{0}^1 \|\nabla F(z_{k+1}+t (z_{k+1}-z_k))\| dt\le \gamma$. It follows from Proposition \ref{prop:minusF} with $\hz=z_{k+1}$ and $z=z_k$ that
\begin{equation}\label{eq:ppm_int}
F(z_{k+1})-F(z_k)=M(z_{k+1}-z_k).
\end{equation}
Therefore, it holds that
\begin{align}
    \begin{split}
         \|\sumF\|^2 & = 2\|F(z_k)\|^2+2\|F(z_{k+1})\|^2-\|\minusF\|^2\\
         &= 2\|F(z_k)\|^2+2\|F(z_{k+1})\|^2-\|M\pran{\minusz}\|^2\\
        &=  2\|F(z_k)\|^2+2\|F(z_{k+1})\|^2-\|sM\Fzko\|^2\\
        &\ge 2\|F(z_k)\|^2+2\|F(z_{k+1})\|^2-\|\Fzko\|^2\\
        &= 2\|F(z_k)\|^2+\|F(z_{k+1})\|^2 \ ,
    \end{split}{}
\end{align}
where the second equality is from the iterate update \eqref{eq:ppm} and the inequality uses $s\le\frac{1}{\gamma}\le \|M\|$. \qed

\vspace{0.2cm}

Let us go back to prove Theorem \ref{thm:ppm}:

{\bf Proof of Theorem \ref{thm:ppm}.}
Let $M=\int_{0}^1 \nabla F(z_k+t (z_{k+1}-z_k))dt$, then it follows from Proposition \ref{prop:minusF} with $\hat{z}=z_{k+1}$ and $z=z_k$ that $F(z_{k+1})-F(z_k)=M(z_{k+1}-z_k)$, thus
\begin{align}\label{eq:Fzko}
    \begin{split}
        F(z_{k+1})&=\frac{1}{2}\left(F(z_k)+F(z_{k+1})\right)+ \frac{1}{2}\pran{\Fzko -\Fzk} \\ 
        & =\frac{1}{2}\left(F(z_k)+F(z_{k+1})\right)+ \frac{1}{2}M\pran{z_{k+1} -z_k} \\
        & = \frac{1}{2}\left(F(z_k)+F(z_{k+1})\right)- \frac{s}{2}M \Fzko \ ,
    \end{split}
\end{align}
where the last equality utilizes the iterate update  \eqref{eq:ppm}. By rearranging \eqref{eq:Fzko}, we obtain
\begin{equation*}
F(z_{k+1})= \frac{1}{2}\pran{I+ \frac{s}{2} M}^{-1}\pran{\sumF} \ ,
\end{equation*}
whereby 
\begin{align}\label{eq:ppm_p0}
    \begin{split}
        \minusF &= M \pran{\minusz} = -sM F(z_{k+1})= -\frac{s}{2}M\pran{I+ \frac{s}{2} M}^{-1}\pran{\sumF} \\
        & = -\frac{s}{2} M \pran{\sum_{i=0}^{\infty} (-1)^i \pran{\frac{s}{2}}^i M^i} \pran{\sumF} \ ,
    \end{split}{}
\end{align}
where the first equality uses \eqref{eq:ppm_int} and the second equality is due to the update rule \eqref{eq:ppm}.

\red{
Going back to the proof scratch stated in Section \ref{sec:scratch-proof}, \eqref{eq:ppm_p0} shows that it holds for PPM that $R(z_k, s)= -\frac{1}{2}M\pran{I+ \frac{s}{2} M}^{-1}$ and $R_i(z_k)= (-1)^{i+1}(\frac{1}{2})^{i+1} M^{i+1}$. The rest of the proof is to show that the $O(s)$-linear-convergence condition \eqref{eq:O(s)-strong-PPM} guarantees the sufficient decay for the corresponding $R_0$ and $R_1$ terms, and the smaller order terms do not affect the rate when the step-size is small enough.}

Notice it holds that
\begin{align}\label{eq:ppm_p1}
    \begin{split}
        &\frac{1}{2}\|\Fzko\|^2-\frac{1}{2}\|\Fzk\|^2 \\
        = & \frac{1}{2}\pran{\sumF}^T \pran{\minusF} \\
        = & -\frac{s}{4} \pran{\sumF}^T  M \sum_{i=0}^{\infty} (-1)^i \pran{\frac{s}{2}}^i M^i \pran{\sumF} \\
        = & \frac{1}{2}\sum_{i=1}^{\infty} (-1)^i \pran{\frac{s}{2}}^{i} \pran{\sumF}^T   M^i \pran{\sumF} \ ,
    \end{split}{}
\end{align}
 where the second equality follows from \eqref{eq:ppm_p0}. 

Since $L(x,y)$ is convex-concave, $M$ is generalized block skew-symmetric. Let us denote $M=\left[\begin{smallmatrix} A & B \\ -B & C \end{smallmatrix}\right]$ and then $M^2=\left[\begin{smallmatrix} A^2- BB^T & AB+BC \\ -B^T A-C B & -B^T B + C^2 \end{smallmatrix}\right]$. It follows Proposition \ref{prop:power-symmetry} that for any power $i$, $M^i$ is also generalized block skew-symmetric, thus the off-diagonal terms cancel out when computing $\pran{\sumF}^T   M^i \pran{\sumF}$. Therefore, it holds that

\begin{equation}\label{eq:ODE_p2}
    \begin{array}{cl}
        & \sum_{i=1}^{2} (-1)^i  \sot^{i-1} \pran{\sumF}^T   M^i \pran{\sumF} \\ \\
        = & -\pran{\sumF}^T \left[\begin{matrix} A- \frac{s}{2} A^2+ \frac{s}{2} BB^T & 0 \\ 0 & C- \frac{s}{2} C^2+ \frac{s}{2} B^T B \end{matrix}\right] \pran{\sumF}\ .        
    \end{array}{}
\end{equation}
Meanwhile, it follows from Proposition \ref{prop:power_bound} with $Q=M$ and $c=s$ that for any $i\ge 3$, 
\begin{align}\label{eq:ODE_p4}
    \begin{split}
        & s^{i-1} |\pran{\sumF}^T M^i \pran{\sumF}| \\
        \le & (i-1) (s\gamma)^{i-2} \pran{\sumF}^T \left[\begin{matrix}  s\gamma A+ s BB^T & 0\\0  & s\gamma C+ s B^T B \end{matrix}\right] \pran{\sumF} \\
        \le & (i-1) (s\gamma)^{i-2} \pran{\sumF}^T \left[\begin{matrix}  A+ s BB^T & 0\\0  & C+ s B^T B \end{matrix}\right] \pran{\sumF} ,
    \end{split}
\end{align}
where the last inequality uses $s\gamma\le 1$. Also notice that $s
\gamma\le \frac{1}{3}$, thus $\sum_{i=3}^{\infty}  \pran{\tfrac{1}{2}}^{i-1}(i-1)(s\gamma)^{i-2} = \frac{1}{2}\pran{ \frac{s\gamma}{2}+ \frac{\frac{s\gamma}{2}}{1-\frac{s\gamma}{2}}}\le \frac{1}{4}$. Therefore, it holds that
\begin{align}\label{eq:ODE_p3}
    \begin{split}
        & \sum_{i=3}^{\infty} (-1)^i  \sot^{i-1} \pran{\sumF}^T M^i \pran{\sumF} \\
        \le & 
        \sum_{i=3}^{\infty}  \pran{\tfrac{1}{2}}^{i-1} s^{i-1} |\pran{\sumF}^T M^i \pran{\sumF}|\\
        \le &  \sum_{i=3}^{\infty}  \pran{\tfrac{1}{2}}^{i-1}(i-1)(s\gamma)^{i-2} \pran{\sumF}^T \conditionmatrix \pran{\sumF} \\
        \le & \frac{1}{4} \pran{\sumF}^T \conditionmatrix \pran{\sumF} \\
        \le & \frac{1}{2} \pran{\sumF}^T \secondconditionmatrix \pran{\sumF} \ ,
    \end{split}{}
\end{align}
where the last inequality follows from $sA^2\preceq s\gamma A \preceq A$ by noticing $A$ is positive semi-definite, $\|A\|\le \|M\| \le \gamma$ and $s\gamma\le 1$.
Substituting \eqref{eq:ODE_p2} and \eqref{eq:ODE_p4} into \eqref{eq:ppm_p1} yields
\begin{equation}\label{eq:ppm_p5}
\begin{array}{cl}
     & \frac{1}{2}\|\Fzko\|^2-\frac{1}{2}\|\Fzk\|^2 \\ \\
     \le & -\frac{s}{8}\pran{\sumF}^T\secondconditionmatrix\pran{\sumF} \\ \\
     \le & -\frac{s\rho(s)}{8} \normt{\sumF} \\ \\
     \le & -\frac{s\rho(s)}{4} \normt{\Fzk} -\frac{s\rho(s)}{8} \|\Fzko\|^2
\end{array}
\end{equation}



where the inequality is due to Proposition \ref{prop:sumF}. By rearranging \eqref{eq:ppm_p5}, we have
$$
\normt{\Fzko} \le \frac{1-\frac{s\rho(s)}{2}}{1+\frac{s\rho(s)}{4}} \normt{\Fzk}\ ,
$$
which furnishes the proof of Theorem \ref{thm:ppm}. \qed




\subsection{Proof of Theorem \ref{thm:slow-mp}}



The following two propositions will be needed for the proof of Theorem \ref{thm:slow-mp}.
\begin{prop}\label{prop:mpFtz}
    Consider EGM with step-size $s$. Let $M=\int_{0}^1 \nabla F(z_{k}+t (z_{k+1}-z_k))dt$, $M_1=\int_{0}^1 \nabla F(\tz_k+t (z_{k+1}-\tz_k))dt$, and $M_2=\int_{0}^1 \nabla F(z_k+t (\tz_{k}-z_k))dt$. Then it holds for any $k$ that
    \begin{equation}\label{eq:mpFztk}
        \Ftz  = \frac{1}{2}\pran{I+ \frac{s}{2} M + \frac{s^3}{2}\threeM}^{-1}\pran{I-\frac{s^2}{2}\twoM}\pran{\sumF} \ .
    \end{equation}
\end{prop}
{\bf Proof.} By the definition of $M$, $M_1$ and $M_2$, we have $\norm{M},\norm{M_1}, \norm{M_2}\le \gamma$. Moreover, it follows from Proposition \ref{prop:minusF} that
\begin{align}
F(z_{k+1})-F(z_k) & =M(z_{k+1}-z_k),\label{eq:mp_zoz}    \\
F(z_{k+1})-F(\tz_k) & =M_1(z_{k+1}-\tz_k), \label{eq:mp_zotz}    \\
F(\tz_{k})-F(z_k) & =M_2(\tz_{k}-z_k) \ , \label{eq:mp_tzz}    
\end{align}

Together with the iterate update of EGM algorithm \eqref{eq:mp}, we obtain
\begin{equation}\label{eq:mp_secondterm}
    \minusF = M \pminusz = -s M\Ftz \ .
\end{equation}
and
\begin{align}\label{eq:mp_thirdterm}
\begin{split}
    \Ftz-\Fzko &= M_1 (\tzk-\zko)= sM_1(\Ftz-\Fzk)=sM_1M_2(\tzk-\zk) =-s^2 M_1M_2\Fzk \\
    & = -s^2 M_1M_2\bracket{\thalf\psumF-\thalf\pminusF} \\
    & = -s^2 M_1M_2\bracket{\thalf\psumF+\thalf s M \Ftz} ,
\end{split}{}
\end{align}
where the second equality is from the update rule \eqref{eq:mp} and the last equality uses \eqref{eq:mp_secondterm}. Using \eqref{eq:mp_secondterm} and \eqref{eq:mp_thirdterm}, we can rewrite $\Ftz$ as:
\begin{align}\label{eq:mp_Ftz_mid}
    \begin{split}
        \Ftz & = \frac{1}{2}\psumF+ \frac{1}{2}\pminusF + \pran{\Ftz-\Fzko} \\
        &= \frac{1}{2}\psumF -\frac{s}{2}M \Ftz - \frac{s^2}{2} M_1M_2 \psumF - \frac{s^3}{2} M_1M_2M \Ftz \ .
    \end{split}
\end{align}
We finish the proof by rearranging \eqref{eq:mp_Ftz_mid}. \qed
\red{
\begin{rem}
Going back to the proof scratch stated in Section \ref{sec:scratch-proof}, Proposition \ref{prop:mpFtz} shows that it holds for EGM that
\begin{align*}
    \minusF &= M\pran{\minusz} = -sM F(z_{k+1})\\
    &=-s\frac{1}{2}M\pran{I+ \frac{s}{2} M + \frac{s^3}{2}\threeM}^{-1}\pran{I-\frac{s^2}{2}\twoM}\pran{\sumF} \ ,
\end{align*}
whereby $R(z_k, s)=-\frac{1}{2}M\pran{I+ \frac{s}{2} M + \frac{s^3}{2}\threeM}^{-1}\pran{I-\frac{s^2}{2}\twoM}$. The rest of the proofs of Theorem \ref{thm:slow-mp} and Theorem \ref{thm:fast-mp} are to show that the $O(s)$-linear-convergence condition \eqref{eq:O(s)-strong-PPM} corresponds to the sufficient decay for the $R_0$ and $R_1$ terms, and the smaller order terms do not affect the rate when the step-size is small enough. Moreover, the difference between the slow rate (Theorem \ref{thm:slow-mp}) and the fast rate (Theorem \ref{thm:fast-mp}) comes from how small the step-sizes need be in order to bound the smaller order terms.
\end{rem}}

\begin{prop}\label{prop:sumF-mp}
Consider EGM with step-size $s$. Suppose $s\le \frac{1}{2\gamma}$, then it holds for any $k$ that
    $$\|\sumF\|^2 \ge \frac{8}{5}\|F(z_k)\|^2+\frac{8}{5}\|F(z_{k+1})\|^2 \ .$$
\end{prop}{}
{\bf Proof.}
It follows from \eqref{eq:mp_zoz} and \eqref{eq:mp} that
\begin{align}\label{eq:mp_temp1}
    \begin{split}
         \|\sumF\|^2 & = 2\|F(z_k)\|^2+2\|F(z_{k+1})\|^2-\|\minusF\|^2\\
         &= 2\|F(z_k)\|^2+2\|F(z_{k+1})\|^2-\|M\pran{\minusz}\|^2\\
        &=  2\|F(z_k)\|^2+2\|F(z_{k+1})\|^2-\|sM\Ftz\|^2 \ .
    \end{split}{}
\end{align}
From Proposition \ref{prop:mpFtz}, we obtain that
\begin{align}\label{eq:mp_temp2}
    \begin{split}
        \normt{sM\Ftz}& \le \tfrac{s^2}{4} \|M\|^2\|I+ \tfrac{s}{2} M + \tfrac{s^3}{2}\threeM\|^{-2}\normt{I-\tfrac{s^2}{2}\twoM}\normt{\sumF}\ \\
        & \le  \tfrac{\pran{\sg}^2}{4} (1-\tfrac{\sg}{2}-\tfrac{\pran{\sg}^3}{2})^{-2} (1+\tfrac{(\sg)^2}{2})^2\normt{\sumF} \\
        & \le \frac{1}{4}\normt{\sumF} \ ,
    \end{split}
\end{align}
where the second inequality comes from the facts: $$\|I+ \tfrac{s}{2} M + \tfrac{s^3}{2}\threeM\| \ge \|I\|-\|\tfrac{s}{2} M\|-\|\tfrac{s^3}{2}\threeM\|
\ge 1-\tfrac{\sg}{2}-\tfrac{\pran{\sg}^3}{2}\ ,$$ and 
$$
\|I-\tfrac{s^2}{2}\twoM\|\le \|I\|+\|\tfrac{s^2}{2}\twoM\|\le 1+\tfrac{(\sg)^2}{2} \ ,
$$
and the last inequality uses the fact that $s\gamma\le \frac{1}{2}$. 
Combining \eqref{eq:mp_temp1} and \eqref{eq:mp_temp2}, we arrive at
$$
\|\sumF\|^2 = 2\|F(z_k)\|^2+2\|F(z_{k+1})\|^2-\|sM\Ftz\|^2 \ge 2\|F(z_k)\|^2+2\|F(z_{k+1})\|^2- \frac{1}{4}\normt{\sumF} \ ,
$$
which finishes the proof by rearrangement. \qed

\vspace{0.2cm}
Let us go back to the proof of Theorem \ref{thm:slow-mp}:

{\bf Proof of Theorem \ref{thm:slow-mp}}
It follows from \eqref{eq:mp_zoz} that
\begin{align}\label{eq:mp_p1}
    \begin{split}
        &\frac{1}{2}\|\Fzko\|^2-\frac{1}{2}\|\Fzk\|^2 \\
        = & \frac{1}{2}\pran{\sumF}^T \pran{\minusF} \\
        = & \frac{1}{2}\pran{\sumF}^T  M \pran{\minusz} \\ 
        = & -\frac{s}{2} \pran{\sumF}^T  M \Ftz \\ 
        = & -\frac{s}{4} \pran{\sumF}^T  M \pran{I+ \frac{s}{2} M + \frac{s^3}{2}\threeM}^{-1}\pran{I-\frac{s^2}{2}\twoM}\pran{\sumF} \\
        = & -\frac{s}{4} \pran{\sumF}^T \pran{ M - \frac{s}{2} M^2} \pran{\sumF} \\
        & -\frac{s}{4} \pran{\sumF}^T \pran{- \frac{s^3}{2}M\threeM} \pran{\sumF} \\
        & -\frac{s}{4} \pran{\sumF}^T  M \sum_{i=2}^{\infty} (-1)^i \pran{\frac{s}{2} M + \frac{s^3}{2}\threeM}^i \pran{\sumF} \\
         & -\frac{s}{4} \pran{\sumF}^T  M \pran{I+ \frac{s}{2} M + \frac{s^3}{2}\threeM}^{-1}\frac{s^2}{2}\twoM\pran{\sumF} \ ,\\
    \end{split}{}
\end{align}
where the third equality is from the update of EGM algorithm, the fourth equality follows from Proposition \ref{prop:mpFtz}, and the last equality is rearrangement by noticing $\pran{I+ \frac{s}{2} M + \frac{s^3}{2}\threeM}^{-1}=\sum_{i=0}^{\infty} (-1)^i \pran{\frac{s}{2} M + \frac{s^3}{2}\threeM}^i$. 

Now let us examine each term at the right-hand-side of \eqref{eq:mp_p1}. In principal, the last three terms are at most $O(s^3)$, and the first term is at least $O(s^2)$, which dominants the right-hand-side of \eqref{eq:mp_p1} when $s$ is small. Suppose $M=\ABC$, then $M^2=\Msquare$. Notice that $\|M_1\|, \|M_2\|, \|M\| \le \gamma \le 1/2s$. For the first term at the right-hand-side of \eqref{eq:mp_p1}, it holds that
\begin{equation}\label{eq:mp_p2}
    \begin{array}{cl}
        & -\frac{s}{4}\pran{\sumF}^T \pran{ M - \frac{s}{2} M^2 } \pran{\sumF} \\ \\
        = & -\frac{s}{4}\pran{\sumF}^T \left[\begin{matrix} A- \frac{s}{2} A^2+ \frac{s}{2} BB^T & 0 \\ 0 & C- \frac{s}{2} C^2+ \frac{s}{2} B^T B \end{matrix}\right] \pran{\sumF} \\    \\  
        \le & -\frac{s\rho(s)}{8} \|\sumF\|^2 \ ,
    \end{array}{}
\end{equation}
where the inequality uses the condition \eqref{eq:O(s)-strong}. For the second term at the right-hand-side of \eqref{eq:mp_p1}, it holds that
\begin{equation}\label{eq:mp_part1}
    \pabs{\frac{s}{4}\pran{\sumF}^T \frac{s^3}{2}M\threeM \pran{\sumF}} \le \frac{s^4}{8} \gamma^4 \|\sumF\|^2 \le \frac{s^3}{16} \gamma^3 \|\sumF\|^2\ ,
\end{equation}
where the last inequality uses $s\gamma\le \frac{1}{2}$. For the third term at the right-hand-side of \eqref{eq:mp_p1}, it holds that
\begin{align}\label{eq:mp_part3}
    \begin{split}
        &\  \pabs{\frac{s}{4}\pran{\sumF}^T  M \sum_{i=2}^{\infty} (-1)^i \pran{\frac{s}{2} M + \frac{s^3}{2}\threeM}^i \pran{\sumF}} \\
        \le &\ \frac{s}{4} \sum_{i=2}^{\infty} \pran{\frac{s}{2} \gamma + \frac{s^3}{2}\gamma^3}^i \gamma \|\sumF\|^2 \\
        \le & \ \frac{s}{4}\sum_{i=2}^{\infty} (\tfrac{5}{8} s\gamma)^i\gamma \|\sumF\|^2 \\
        = & \ \frac{25}{256} s^3\gamma^3 \frac{1}{1-\frac{5}{8} s\gamma} \|\sumF\|^2 \\
        \le & \frac{5}{32} s^3\gamma^3 \|\sumF\|^2 \ ,
    \end{split}
\end{align}
where the first inequality is because $$\norm{\sum_{i=2}^{\infty} (-1)^i \pran{\frac{s}{2} M + \frac{s^3}{2}\threeM}^i} \le \sum_{i=2}^{\infty}  \pran{\frac{s}{2} \norm{M} + \frac{s^3}{2}\norm{\threeM} }^i=\sum_{i=2}^{\infty} \pran{\frac{s}{2} \gamma + \frac{s^3}{2}\gamma^3}^i\ ,$$
and the second and last inequality uses the fact that $s\gamma\le \tfrac{1}{2}$. Similarly, for the last term at the right-hand-side of \eqref{eq:mp_p1}, it holds that
\begin{align}\label{eq:mp_part2}
    \begin{split}
        &\  \pabs{\frac{s^3}{8} \pran{\sumF}^T  M \pran{I+ \frac{s}{2} M + \frac{s^3}{2}\threeM}^{-1}\twoM\pran{\sumF}} \\
        \le & \ \frac{s^3\gamma^3}{8} \frac{1}{1-\frac{s\gamma}{2}-\frac{s^3\gamma^3}{2}}\|\sumF\|^2 \\
        \le & \ \frac{1}{5}s^3\gamma^3\|\sumF\|^2 \ .
    \end{split}{}
\end{align}

Substituting \eqref{eq:mp_p2}, \eqref{eq:mp_part1}, \eqref{eq:mp_part2} and \eqref{eq:mp_part3} into \eqref{eq:mp_p1}, we arrive at:
\begin{align}\label{eq:mp-almostdone}
    \begin{split}
        & \frac{1}{2}\|\Fzko\|^2-\frac{1}{2}\|\Fzk\|^2\\
        \le & \pran{-\frac{s\rho(s)}{8} +\pran{\frac{1}{16}+\frac{5}{32}+\frac{1}{5}}s^3\gamma^3}\normt{\sumF} \\
        \le & \pran{-\frac{s\rho(s)}{8} +\frac{1}{2} s^3\gamma^3}\normt{\sumF} \\
        \le & -\frac{s\rho(s)}{16} \normt{\sumF} \\
        \le & -\frac{s\rho(s)}{10} \|\Fzko\|^2 -\frac{s\rho(s)}{10}\|\Fzk\|^2\ ,
    \end{split}
\end{align}
where the third inequality uses $\rho(s)\ge 8s^2\gamma^3$, and the last inequality is from Proposition \ref{prop:sumF-mp}. Rearranging \eqref{eq:mp-almostdone} yields
\begin{equation*}
    \normt{\Fzko} \le \pran{\frac{1-\frac{s\rho(s)}{5}}{1+\frac{s\rho(s)}{5}}} \normt{\Fzk}\ ,
\end{equation*}
which finishes the proof by telescoping.\qed


\subsection{Proof of Theorem \ref{thm:fast-mp}}
\vspace{0.2cm}
The next proposition will be used in the proof of Theorem \ref{thm:fast-mp}:
\begin{prop}\label{prop:cQc}
Consider $Q\in\RR^{(m+n)\times(m+n)}$ with $\|Q\|\le \alpha < 1$. Suppose there exist a positive semi-definite matrix $P$ satisfies that for any $c\in \RRmn$ and any positive integer $k\ge 3$, it holds that 
\begin{equation}\label{eq:prop_condition}
    |c^T Q^k c|\le (k-1)\alpha^{k-2}s^2 c^T P c
\end{equation}
with a positive scalar $s$, then we have for any $j\ge 3$ that
\begin{equation}
    \left|c^T Q^j(I+\frac{Q}{2}+\frac{Q^3}{2})^{-1}c\right|\le s^2 h_2(2\alpha)(2\alpha)^{j-2} {c^TPc} \ ,
\end{equation}
where $h_2(u)=\pran{1-\frac{u}{2}-\frac{u^3}{2}}^{-1}$.
\end{prop}
\textbf{Proof.}
Consider function $h_1(u):=(1+\frac{u}{2}+\frac{u^3}{2})^{-1}$ and $h_2(u):=(1-\frac{u}{2}-\frac{u^3}{2})^{-1}$. The power series expansion of $h_1(u)$ and $h_2(u)$ are
\begin{equation}\label{eq:sub_condition}
    h_1(u)=\pran{1+\frac{u}{2}+\frac{u^3}{2}}^{-1}= \sum_{l=0}^{\infty}(-1)^{l} \pran{\frac{u}{2}+\frac{u^3}{2}}^l=\suminf a_i u^i \ ,
\end{equation}
and 
\begin{equation}{\label{eq:g_expansion}}
    h_2(u)=\pran{1-\frac{u}{2}-\frac{u^3}{2}}^{-1}= \sum_{l=0}^{\infty} \pran{\frac{u}{2}+\frac{u^3}{2}}^l=\suminf b_i u^i \ ,
\end{equation}
where $a_i$ and $b_i$ are the $i$-th coefficients of the power series expansion of $h_1(u)$ and $h_2(u)$, respectively. Notice that the above two infinite sum converges in the domain $\{u:|\frac{u}{2}+\frac{u^3}{2}|<1\}$. Furthermore, it is straight-forward to see that for any $i$, $|a_i|\le b_i$ because of the existence of the $(-1)^l$ term in the expansion of $h_1(u)$.




Notice that $\|Q\|\le \alpha < 1$, thus $\|\frac{Q}{2}+\frac{Q^3}{2}\|<1$, whereby the power series expansion of the matrix function $f(Q)$ converge. Therefore, it holds that
\begin{align}\label{eq:temp1}
    \begin{split}
        \left|c^T Q^j\pran{I+\frac{Q}{2}+\frac{Q^3}{2}}^{-1}c\right|&=\pabs{c^T \suminf a_i Q^{i+j}c}\le \suminf \pabs{a_i}\pabs{c^TQ^{i+j}c} \le \suminf \pabs{a_i}(i+j-1)\alpha^{i+j-2}s^2{c^TPc} \ , \\
    \end{split}
\end{align}
where the last inequality is from \eqref{eq:prop_condition}. Furthermore, notice that $j\ge 3$, thus it holds for any $i\ge 0$ that $(i+j-1)\alpha^{i+j-2}\le (2\alpha)^{i+j-2}$. Therefore,
\begin{align}\label{eq:temp2}
    \begin{split}
        \suminf \pabs{a_i}(i+j-1)\alpha^{i+j-2}{c^TPc}\le & \suminf \pabs{a_i}(2\alpha)^{i+j-2}{c^TPc}\le \suminf b_i(2\alpha)^{i+j-2}{c^TPc}
        =   h_2(2\alpha)(2\alpha)^{j-2} {c^TPc}  \ ,
    \end{split}
\end{align}
where the second inequality uses $|a_i|\le b_i$, the first equality is from \eqref{eq:g_expansion}. Combining \eqref{eq:temp1} and \eqref{eq:temp2} finishes the proof of Proposition \ref{prop:cQc}. \qed

Now let us go back to EGM. By choosing $Q=sM$, $\alpha=s\gamma$, and $P=\middlematrix$ in Proposition \ref{prop:cQc}, we obtain:
\begin{cor}\label{cor:inverse}
\begin{equation}\label{eq:inverse}
     \pabs{s^jc^TM^j\pran{I+ \frac{s}{2} M + \frac{s^3}{2}M^3}^{-1}c} \le s^2 (1-s\gamma-4s^3\gamma^3)^{-1}(2s\gamma)^{j-2}{c^T \middlematrix c} \ .
\end{equation}
\end{cor}
\textbf{Proof.} Notice that $\|sM\|\le s\gamma < 1$. Furthermore, it follows by Proposition \ref{prop:power_bound} that for any $c$ and $k\ge 3$,
$$
|c^T s^k M^k c|= s^k |c^T M^k c|\le (k-1)s^2(s\gamma)^{k-2}c^T\middlematrix c \ .
$$
Thus $Q=sM$, $\alpha=s\gamma$, and $P=\middlematrix$ satisfies the conditions in Proposition \ref{prop:cQc}, which leads to \eqref{eq:inverse} by noticing $h_2(2s\gamma)=(1-s\gamma-4s^3\gamma^3)^{-1}$. \qed

\vspace{0.2cm}
\textbf{Proof of Theorem \ref{thm:fast-mp}.}
Following the notations in the proof of Theorem \ref{thm:slow-mp}, it holds that $M_1=M_2=M=\ABC$ when the minimax function $L(x,y)$ is quadratic, and we can then write \eqref{eq:mp_p1} as
\begin{align}\label{eq:mp_f_p0}
    \begin{split}
        &\frac{1}{2}\|\Fzko\|^2-\frac{1}{2}\|\Fzk\|^2 \\
        = & -\frac{s}{4} \pran{\sumF}^T \pran{ M - \frac{s}{2} M^2} \pran{\sumF} \\
        & +\frac{s^4}{8} \pran{\sumF}^T M^4 \pran{\sumF} \\
        & -\frac{s}{4} \pran{\sumF}^T  M \sum_{i=2}^{\infty} (-1)^i \pran{\frac{s}{2} M + \frac{s^3}{2}M^3}^i \pran{\sumF} \\
         & -\frac{s^3}{8} \pran{\sumF}^T  M^3 \pran{I+ \frac{s}{2} M + \frac{s^3}{2}M^3}^{-1}\pran{\sumF} \ ,\\
    \end{split}{}
\end{align}
by utilizing the fact that $f(M) M=M f(M)$ if $f$ is a function of $M$  with convergent power series. Let us again examine each term at the right-hand side of \eqref{eq:mp_f_p0}. For the first term, recall that \eqref{eq:mp_p2} shows that 
\begin{align}\label{eq:mp_f_p1}
    \begin{split}
        &  -\frac{s}{4}\pran{\sumF}^T \pran{ M - \frac{s}{2} M^2 } \pran{\sumF} \\
        = & -\frac{s}{4}\pran{\sumF}^T \secondconditionmatrix \pran{\sumF} \ .
    \end{split}{}
\end{align}
For the second term, it follows from Proposition \ref{prop:power_bound} that
\begin{align}\label{eq:mp_f_p2}
    \begin{split}
        &\ \frac{s^4}{8} \pabs{\pran{\sumF}^T M^4 \pran{\sumF}} \\ \le & \  \frac{3s^4}{8}\gamma^2 {\pran{\sumF}^T \middlematrix \pran{\sumF}} \\
        \le & \  \frac{3s}{8}(s\gamma)^2 {\pran{\sumF}^T \left[\begin{matrix}   A+ s BB^T & 0\\0  &  C+ s B^T B \end{matrix}\right] \pran{\sumF}}\\
        \le & \  \frac{3s}{4}(s\gamma)^2 {\pran{\sumF}^T \secondconditionmatrix \pran{\sumF}} \ .
    \end{split}
\end{align}
For the third term, it holds that
\begin{align}\label{eq:mp_f_p3}
\begin{split}
        &~\pabs{\frac{s}{4} \pran{\sumF}^T  M \sum_{i=2}^{\infty} (-1)^i \pran{\frac{s}{2} M + \frac{s^3}{2}M^3}^i \pran{\sumF}} \\
    = & ~\pabs{\frac{s}{4} \pran{\sumF}^T  M \pran{\frac{s}{2} M + \frac{s^3}{2}M^3}^2\sum_{i=0}^{\infty} (-1)^i \pran{\frac{s}{2} M + \frac{s^3}{2}M^3}^i \pran{\sumF}} \\
    =& ~\pabs{\frac{s}{4} \pran{\sumF}^T  M \pran{\frac{s}{2} M + \frac{s^3}{2}M^3}^2 \pran{I+ \frac{s}{2} M + \frac{s^3}{2}M^3}^{-1} \pran{\sumF}} \\
    =&~ \pabs{\frac{s}{4} \pran{\sumF}^T  M \pran{\frac{s^2}{4}M^2+\frac{s^4}{2}M^4+\frac{s^6}{4}M^6} \pran{I+ \frac{s}{2} M + \frac{s^3}{2}M^3}^{-1} \pran{\sumF}} \\
    \le &~ \tfrac{s^2}{4}    \pran{\tfrac{(2s\gamma)}{4}+\tfrac{(2s\gamma)^3}{2}+\tfrac{(2s\gamma)^5}{4}} \pran{1 -s\gamma - 4s^3\gamma^3}^{-1} \times \\
    &~ \pran{\sumF}^T \middlematrix \pran{\sumF} \\
    \le &~ \tfrac{s}{4}    \pran{\tfrac{(2s\gamma)}{4}+\tfrac{(2s\gamma)^3}{2}+\tfrac{(2s\gamma)^5}{4}} \pran{1 -s\gamma - 4s^3\gamma^3}^{-1} \times \\
    &~ \pran{\sumF}^T \conditionmatrix \pran{\sumF} \\
    \le &~ \tfrac{s}{2}    \pran{\tfrac{(2s\gamma)}{4}+\tfrac{(2s\gamma)^3}{2}+\tfrac{(2s\gamma)^5}{4}} \pran{1 -s\gamma - 4s^3\gamma^3}^{-1} \times \\
    &~ \pran{\sumF}^T \secondconditionmatrix \pran{\sumF} \ ,
\end{split}
\end{align}
where the second equality is because $\pran{I+ \frac{s}{2} M + \frac{s^3}{2}M^3}^{-1}=\sum_{i=0}^{\infty} (-1)^i \pran{\frac{s}{2} M + \frac{s^3}{2}M^3}^i$, the first inequality utilizes Corollary  \ref{cor:inverse}, the second inequality uses $s\gamma \le 1$. 

For the fourth term, it follows Corollary  \ref{cor:inverse} that

\begin{align}\label{eq:mp_f_p4}
    \begin{split}
        &~\pabs{\frac{s^3}{8} \pran{\sumF}^T  M^3 \pran{I+ \frac{s}{2} M + \frac{s^3}{2}M^3}^{-1}\pran{\sumF}} \\
        \le &~ \frac{s^2}{8} (2s\gamma) \pran{1 -s\gamma - 4s^3\gamma^3}^{-1} \pran{\sumF}^T \middlematrix \pran{\sumF}\\
         \le &~ \frac{s}{8} (2s\gamma) \pran{1 -s\gamma - 4s^3\gamma^3}^{-1} \pran{\sumF}^T \conditionmatrix \pran{\sumF} \\
         \le &~ \frac{s}{4} (2s\gamma) \pran{1 -s\gamma - 4s^3\gamma^3}^{-1} \pran{\sumF}^T \secondconditionmatrix \pran{\sumF} \ .
    \end{split}
\end{align}
Substituting \eqref{eq:mp_f_p1}, \eqref{eq:mp_f_p2}, \eqref{eq:mp_f_p3}, \eqref{eq:mp_f_p4} into \eqref{eq:mp_f_p0}, and noticing that $s\gamma\le \frac{1}{8}$, we obtain
\begin{align}\label{eq:mp_f_p5}
    \begin{split}
        &\frac{1}{2}\|\Fzko\|^2-\frac{1}{2}\|\Fzk\|^2 \\
        \le &- \frac{s}{4} \pran{1-3(s\gamma)^2-2\pran{\tfrac{(2s\gamma)}{4}+\tfrac{(2s\gamma)^3}{2}+\tfrac{(2s\gamma)^5}{4}} \pran{1 -s\gamma - 4s^3\gamma^3}^{-1} - (2s\gamma) \pran{1 -s\gamma - 4s^3\gamma^3}^{-1}}\times \\
        & \pran{\sumF}^T \secondconditionmatrix \pran{\sumF} \\
        \le &  -\frac{s}{8} \pran{\sumF}^T \secondconditionmatrix \pran{\sumF} \\
        \le & -\frac{s\rho(s)}{16} \normt{\sumF} \ .
    \end{split}{}
\end{align}
It then follows from Proposition \ref{prop:sumF-mp} that
$$
\frac{1}{2}\|\Fzko\|^2-\frac{1}{2}\|\Fzk\|^2 \le -\frac{s\rho(s)}{10} \|\Fzko\|^2 -\frac{s\rho(s)}{10}\|\Fzk\|^2 \ ,
$$
and after rearrangement, we arrive at 
\begin{equation*}
    \normt{\Fzko} \le \pran{\frac{1-\frac{s\rho(s)}{5}}{1+\frac{s\rho(s)}{5}}} \normt{\Fzk}\ ,
\end{equation*}
which finishes the proof by telescoping. \qed

\end{document}